\documentclass[preprint,3p,times]{elsarticle}

\usepackage{amssymb}
\usepackage{amsthm}
\usepackage{graphicx}
\usepackage[utf8]{inputenc}
\usepackage{amsmath}
\usepackage{amssymb}
\usepackage{mathtools}
\usepackage{bbm}
\usepackage{cleveref}
\usepackage[font=small,labelfont={normal, bf}]{caption}
\usepackage{subcaption}
\usepackage{comment}
\usepackage{svg}
\makeatletter
\newcommand{\@giventhatstar}[2]{#1\;\middle|\;#2)}
\newcommand{\@giventhatnostar}[3][]{#1#2\;#1|\;#3#1}
\newcommand{\giventhat}{\@ifstar\@giventhatstar\@giventhatnostar}
\makeatother
\newtheorem{theorem}{Theorem}[section]
\newtheorem{fact}{Fact}[section]
\theoremstyle{definition}
\newtheorem{definition}{Definition}[section]
\newtheorem{corollary}{Corollary}[section]
\newtheorem{example}{Example}[section]

\biboptions{sort&compress}

\begin{document}

\begin{frontmatter}

\title{Scaled Brownian motion with random anomalous diffusion exponent}

\author[inst1]{Hubert Woszczek}
\ead{hubert.woszczek@pwr.edu.pl}
\author[inst1,inst2,inst3]{Aleksei Chechkin}
\author[inst1]{Agnieszka Wyłomańska}

\affiliation[inst1]{organization={Faculty of Pure and Applied Mathematics, Hugo Steinhaus Cenetr, Wroclaw University of Science and Technology},
            addressline={Wyspiańskiego 27}, 
            city={Wrocław},
            postcode={50-370}, 
            country={Poland}}

\affiliation[inst2]{organization={Institute for Physics \& Astronomy, University of Potsdam},
            addressline={}, 
            city={Potsdam-Golm},
            postcode={14476 }, 
            country={Germany}}

\affiliation[inst3]{organization={Akhiezer Institute for Theoretical Physics, National Science Center ‘Kharkov Institute of Physics and Technology’,},
            addressline={Akademicheskaya st.1}, 
            city={Kharkov},
            postcode={61108}, 
            country={Ukraine}}


\begin{abstract}
The scaled Brownian motion (SBM) is regarded as one of the paradigmatic random processes, featuring the anomalous diffusion property characterized by the  diffusion exponent. It is a Gaussian, self-similar process with independent increments, which has found applications across various fields, from turbulence and stochastic hydrology to biophysics. In our paper, inspired by recent single particle tracking biological experiments, we introduce a process termed the scaled Brownian motion with random exponent (SBMRE), which preserves SBM characteristics at the level of individual trajectories, albeit with randomly varying anomalous diffusion exponents across the trajectories.
We discuss the main probabilistic properties of SBMRE, including its probability density function (pdf), and the q-th absolute moment. Additionally, we present the expected value of the time-averaged mean squared displacement (TAMSD) and the ergodicity breaking parameter. Furthermore, we analyze the pdf of the first hitting time in a semi-infinite domain, the martingale property of SBMRE, and its stochastic exponential. As special cases, we consider two distributions of the anomalous diffusion exponent, namely the two-point and beta distributions, and discuss the asymptotics of the presented characteristics in such cases. Theoretical results for SBMRE are validated through numerical simulations and compared with the corresponding characteristics for SBM.

\end{abstract}

\end{frontmatter}

\section{Introduction}

The stochastic processes of anomalous diffusion behavior are ubiquitous in many areas of interest. They are characterized  by the non-linear, mostly power-law, second moment, namely $ \mathbb{E}\left[ X^2(t)\right] \sim t^\mu$ \cite{Bouchaud1990}. The $\mu$ parameter is called the anomalous diffusion exponent. Depending on the $\mu$ value, one can distinguish between sub-diffusive ($\mu
< 1$) and super-diffusive ($\mu > 1$) behavior. When $\mu=1$, the process exhibits linear time dependence that is the hallmark of the normal (ordinary) Brownian diffusion \cite{Mazo2008}.

Mathematically, ordinary Brownian motion, or Wiener process, is the Gaussian self-similar process with stationary independent increments \cite{Schilling2012}. There are two generic Gaussian processes exhibiting anomalous diffusion behavior. A classic example is the fractional Brownian motion (FBM), introduced by Kolmogorov \cite{kol140} and further examined by Mandelbrot and van Ness \cite{NessMandelbrot}, particularly within the context of economic time series analysis. FBM is characterized by the Hurst exponent $H\in(0,1)$, which dictates its anomalous diffusion behavior as $\mu=2H$. It is the only Gaussian self-similar process with stationary power-law correlated increments, closely associated with what is known as long-range phenomena \cite{beran}. FBM has found numerous interesting applications across various areas, including hydrology \cite{https://doi.org/10.1029/97WR01982,BENSON2013479}, telecommunications and signal processing \cite{8246429,7448970,984735}, image analysis \cite{605414,6879494}, economics \cite{ROSTEK201330,10.1063/5.0054119,XIAO2010935}, biological systems such as single-particle tracking experiments \cite{weiss2012,franosch2013,metzler2014,szarek2022statistical,gleb2019}, and many others.

Another Gaussian process demonstrating anomalous diffusion behavior is the scaled Brownian motion (SBM) that was introduced implicitly in the context of turbulent diffusion in \cite{Batchelor1952}. It was then rediscovered in \cite{Lim2002} as a simple time-changed subordination of the Brownian motion by a deterministic power function. This process has similar characteristics as FBM, namely it is Gaussian and self-similar, however, in contrast to FBM its increments are independent (which implies Markovianity) and non-stationary.

Many authors have examined the SBM from a theoretical perspective. For instance, in \cite{Debicki1998,Debicki2001}, discussions revolve around the asymptotics of the supremum of SBM, while \cite{Magdziarz2020} delves into the Lamperti transformation for SBM. In \cite{Thiel2014} SBM is considered as the mean-field model for continuous-time random walk. In \cite{Jeon2014} highly nonstationary behavior of the SBM in confining potential is demonstrated, while the authors of \cite{Jeon2014, Thiel2014, Safdari2015} analyze ergodicity breaking in SBM. In \cite{Sposini2019} the power spectral density of SBM is explored. The literature also addresses various modifications of the SBM. For instance, the papers \cite{Bodrova2015,Bodrova2016} explore ultraslow and underdamped versions of SBM, respectively, while \cite{Bodrova2019,Bodrova2019a} discuss the SBM with resetting. The authors of \cite{Lee2019} propose a modification of SBM to render it stationary. Further studies include SBM with random diffusion coefficient \cite{Li2024}, SBM in a quenched
disorder environment \cite{Suleiman2024},  and scaled geometric Brownian motion \cite{Cherstvy2021}. Recently, the authors of \cite{ValdesGomez2023} explore FBM and SBM on a sphere, investigating the effects of long-time correlations on navigation strategies. The SBM has also found intriguing applications, particularly in physical sciences and biology, as evidenced by works such as \cite{Jeon2014,Molini2011,Bassler2007,Bodrova2015,Novikov2014}. Interested readers can also explore \cite{Thapa2022,Grzesiek2018,Balcerek2021} for different approaches proposed for testing SBM.

With the recent developments in single particle tracking technique it became clear that the "standard" anomalous diffusion models, like e.g., FBM and SBM, may be inadequate for describing certain complex systems. For example, modern experiments suggest that the motion of biological cells displays anomalous diffusion behavior at the individual trajectory level. However, there is variability in the anomalous diffusion exponent from one trajectory to another, see e.g.  \cite{wang2018,benelli2021sub,speckner2021single,cherstvy2019non}. In such a scenario the natural generalization of the standard anomalous diffusion model is to take into account the randomness of the parameters characterizing the random process, such as diffusivity and/or anomalous diffusion exponent. The idea of such doubly stochastic behavior is the essence of the superstatistics approach \cite{s_47,Beck_2005} and recently has got its further development in the concept of "diffusing diffusivity" suggested in \cite{Chubynsky2014} and then explored in a set of papers, see e.g., the review \cite{West2023} and references therein.

The multifractional generalization of FBM, which allows the Hurst exponent to be a stationary random process was discussed in mathematical literature \cite{levy1995multifractional,ayachetaqqu2005}. Its simplified version, called FBMRE, which considers the Hurst exponent as a variable randomly changing from trajectory to trajectory is analysed in \cite{Balcerek2022}, where interesting effects like accelerating diffusion and persistence transitions in course of time is found. We also refer the reader to recent publications that discuss FBM-based models with varying scenarios for random parameters \cite{han2020deciphering,korabel2021local,balcerek2023modelling,wang2023memorymultifractional}.

In this paper, we investigate scaled Brownian motion with a random exponent (SBMRE), which was introduced in \cite{Santos2022}. We extend the mathematical description of SBMRE for general distribution of the anomalous diffusion exponent. We present the probability density function (pdf) and analyse q-th  absolute moment of SBMRE, expected value of the time-averaged mean squared displacement (TAMSD), ergodicity breaking parameter, and the first hitting time pdf. Finally, we also discuss martingale property of SBMRE and its stochastic exponential. 
As a special case we examine two generic distributions of the anomalous diffusion exponent, namely two-point and beta distributions. In such cases we present the abovementioned characteristics of SBMRE and analyze their asymptotics. Our theoretical results are supported by numerical analysis demonstrating the specific behavior of the considered stochastic process and its differences in contrast to the classical SBM.

The rest of the paper is organized as follows. In Section \ref{preliminar} we recall all the definitions needed. In Section \ref{sbmsec} we remind the readers important properties of SBM while in Section \ref{sbmresec} we present the same characteristics for SBMRE for general distribution of the anomalous diffusive exponent. In Section \ref{examplessec} we present the results for two exemplary distributions of the anomalous diffusive exponent. In Section \ref{nasec} we confirm the theoretical results with numerical simulations. Last section concludes the paper.

\section{Preliminaries}\label{preliminar}
\noindent In this section, we introduce all  definitions and notations used in the following sections. We assume that all processes take values in $\mathbb{R}$. If we do not express it explicitly, we assume, that corresponding probability space is $\left( \Omega, \mathcal{F}, \mathbb{P}, \left\{\mathcal{F}_{t}\right\}\right)$.

\begin{definition} 
    Let $\left\{ X\left(t\right)\right\}_{t\geq 0}$ be the stochastic process . Then, the time averaged mean square displacement (TAMSD) of the process $\left\{ X\left(t\right)\right\}_{t\geq 0}$ is defined as follows \cite{metzler2014}:
    \begin{equation}\label{tamsd}
    \delta\left(\tau\right) = \frac{1}{T - \tau} \int_0^{T-\tau} \left(X\left(t + \tau\right) - X\left(t\right)\right)^2 dt,
\end{equation}
where $\tau \in \left[0, T\right)$ is the time lag of measured time series (the width of a sliding window), and $T \in \left(0, \infty\right)$ is a time horizon (the trajectory length). 
\end{definition}
\noindent Next, we introduce an important concept, called an ergodicity breaking parameter, which helps to quantify the ergodicity of a given process.
\begin{definition}
Let $\left\{ X\left(t\right)\right\}_{t\geq 0}$ be the stochastic process. Then, the ergodicity breaking parameter of the process $\left\{X\left(t\right)\right\}_{t\geq 0}$ is defined as follows \cite{Rytov1987, He2008}:
\begin{equation}\label{EB_general}
    EB\left(\tau\right) = \frac{\mathbb{E}\left[\delta^2 \left(\tau\right)\right] - \mathbb{E}\left[\delta \left(\tau\right)\right]^2}{\mathbb{E}\left[\delta \left(\tau\right)\right]^2} = \frac{Var\left(\delta\left(\tau\right)\right)}{\mathbb{E}\left[\delta\left(\tau\right)\right]^2} = \frac{\mathcal{N}\left(\tau\right)}{\mathcal{D}\left(\tau\right)}.
\end{equation}
\end{definition}
\noindent Next, we define the hitting time in the barrier $b>0$ of a given stochastic process. 
\begin{definition}
    Let $\left\{ X\left(t\right)\right\}_{t\geq 0}$ be the stochastic process. Then the hitting time in the barrier $b>0$ of the process $\left\{ X\left(t\right)\right\}_{t\geq 0}$ is a random variable defined as follows:
    \begin{equation}
        \tau_b = \inf\left\{t \geq 0: X\left(t\right) = b\right\}.
    \end{equation}
\end{definition}
\noindent The next two definitions are the quadratic variation and the stochastic exponent. These are important tools when we focus on the martingale properties of the process.
\begin{definition}
    Let $\left\{ X\left(t\right)\right\}_{t\geq 0}$ be the semimartingale (for more details, see \cite{Oeksendal2003}). Then the quadratic variation of the process $\left\{X\left(t\right)\right\}_{t \geq 0}$ is defined by \cite{Oeksendal2003}:
    \begin{equation}
         \left[X\right]_t = \lim_{\Delta t_k \to 0} \sum_{t_k \leq t} \left(X\left(t_{k+1}\right) - X\left(t_k\right)\right)^2,
    \end{equation}
    where $0=t_1<t_2< \ldots < t_n = t$, $\Delta t_k = t_{k+1} - t_k$ and limit is taken in probability.
\begin{definition}
    Let $\left\{ X\left(t\right)\right\}_{t\geq 0}$ be the semimatringale. Then the stochastic exponential of the process $\left\{X\left(t\right)\right\}_{t \geq 0}$ is given by \cite{Oeksendal2003}: 
    \begin{equation}
        Y\left(t\right) =\exp\left\{\lambda X\left(t\right) - \frac{\lambda^2}{2} \left [ X\right]_t\right\}, \; \lambda>0.
    \end{equation}
\end{definition}
\end{definition}
\noindent The last two definitions are the big $\mathcal{O}$ and the little $o$ notation. They are important for providing asymptotic formulas. In both definitions, we consider functions $f\left(\cdot\right)$ and $g\left(\cdot\right)$ whose domain and co-domain are $\mathbb{R}$.
\begin{definition}
    We say, that $f\left(x\right) \in \mathcal{O}\left(g\left(x\right)\right)$ as $x \to \infty$, if there exist $c>0$ and $x_0 \in \mathbb{R}$, such that for all $x\geq x_0$ we have \cite{Bruijn1981}: 
    \begin{equation}
        \left|f\left(x\right)\right| \leq c \cdot g\left(x\right).
    \end{equation}
\end{definition}
\begin{definition}
    We say, that $f\left(x\right) \in o\left(g\left(x\right)\right)$ as $x \to \infty$, if for all $c>0$ there exists $x_0 \in \mathbb{R}$, such that for all $x\geq x_0$ we have \cite{Bruijn1981}: 
    \begin{equation}
        \left|f\left(x\right)\right| \leq c \cdot g\left(x\right).
    \end{equation}
\end{definition}
\section{Scaled Brownian motion}\label{sbmsec}
\noindent In this section we recall the definition of the scaled Brownian motion and describe its main properties.
\begin{definition}
Let $\left\{B\left(t\right)\right\}_{t\geq0}$ be the standard Brownian motion. Then, the scaled Brownian motion (SBM) is defined as follows:
\begin{equation}
    B_{\alpha}\left(t\right) = B\left(t^{\alpha}\right),
\end{equation}
where $\alpha>0$ is an anomalous diffusion exponent.
It is a Gaussian self-similar process with independent increments (Markovianity), zero mean, and autocovariance $Cov\left(B_{\alpha}\left(s\right), B_{\alpha}\left(t\right)\right) = \min\left\{t^{\alpha}, s^{\alpha}\right\}$. SBM is subdiffusive process for $\alpha<1$ and superdiffusive process for $\alpha>1$. For $\alpha=1$ the SBM becomes an ordinary Brownian motion.
The probability density function is given by:
\begin{equation}
    p\left(x, t \right) = \frac{1}{\sqrt{2\pi t^{\alpha}}} \exp\left\{-\frac{x^2}{2t^{\alpha}}\right\}, \quad x \in \mathbb{R}.
\end{equation}
The SBM can also be represented by the following stochastic differential equation (SDE)  \cite{Lim2002}:
\begin{equation}\label{sdesbm}
    dB_{\alpha}\left(t\right) = \sqrt{\alpha} t^{\frac{\alpha - 1}{2}}dB\left(t\right), \quad B_{\alpha}\left(0\right)=0.
\end{equation}
The solution of \eqref{sdesbm} reads
\begin{equation}
    B_{\alpha}\left(t\right) = \int_0^t \sqrt{\alpha} s^{\frac{\alpha - 1}{2}}dB\left(s\right).
\end{equation}
\end{definition}
\noindent In the next part of this section we present the main properties of SBM.
\begin{fact}
    Let $\left\{B_{\alpha}\left(t\right)\right\}_{t\geq0}$ be the SBM. Then, the expected value of the TAMSD for SBM is given by \cite{Jeon2014}:
    \begin{equation}\label{tamsdsbm}
        \mathbb{E}\left[\delta_{\alpha}\left(\tau\right)\right] = \frac{T^{\alpha +1} - \tau^{\alpha + 1} - \left(T - \tau\right)^{\alpha + 1}}{\left(\alpha+1\right)\left(T-\tau\right)}.
    \end{equation}
\end{fact}
\noindent For $\tau/T<<1$ the expected value of the TAMSD for SBM \eqref{tamsdsbm} asymptotically behaves like \cite{Safdari2015}: 
\begin{equation}\label{tamsdsbmasympt}
\mathbb{E}\left[\delta_{\alpha}\left(\tau\right)\right] \sim \tau/T^{1-\alpha}.
\end{equation}
\begin{fact}
    Let $\left\{B_{\alpha}\left(t\right)\right\}_{t\geq0}$ be the SBM. Then, the ergodicity breaking parameter for SBM is given by \cite{Safdari2015}
    \begin{equation}\label{ebsbm}
        EB_{\alpha}\left(\tau\right) = \frac{\mathcal{N}_{\alpha}\left(\tau\right)}{\mathcal{D}_{\alpha}\left(\tau\right)},
    \end{equation}
    where 
    \begin{equation}\label{sbmnumerator}
    \begin{aligned}
        \mathcal{N}_{\alpha}\left(\tau\right) = \frac{4 \tau^{2\alpha + 2}}{\left(T- \tau\right)^2}\left[\frac{\left(T/\tau - 1\right)^{2\alpha + 1}}{2\alpha + 1} + \frac{\left(3\alpha+1\right)\left(T/\tau - 1\right)^{2\alpha + 2}}{2\left(\alpha + 1\right)^2 \left(2\alpha + 1\right)} - \frac{2\left(T/\tau\right)^{\alpha + 1}\left(T/\tau - 1\right)^{\alpha + 1}}{\left(\alpha + 1\right)^2}\right. + \\
        + \left.\frac{\left(T/\tau\right)^{2\alpha + 2}}{2\left(\alpha + 1\right) \left(2\alpha + 1\right)} - \frac{\left(2\alpha^2 + \alpha+1\right)}{2\left(\alpha + 1\right)^2 \left(2\alpha + 1\right)} + \frac{2}{\alpha +1} \int_0^{T/\tau - 1} x^{\alpha + 1}\left(x+1\right)^{\alpha}dx\right]
    \end{aligned}
    \end{equation}
    and $\mathcal{D}_{\alpha}\left(\tau\right)$ is the square of \eqref{tamsdsbm}.
\end{fact}
\noindent For $\tau/T<<1$ the EB paramater for SBM \eqref{ebsbm} asymptotically behaves like \cite{Safdari2015}: 
\begin{equation}\label{ebsbmsympt}
    EB_{\alpha}\left(\tau\right)\sim
\begin{cases}
C\left(\alpha\right)\left(\frac{\tau}{T}\right)^{2\alpha}; \quad 0<\alpha<1/2\\
\frac{\tau}{12T}\left[\log\left(\frac{T}{\tau}\right) + 2\log\left(2\right) - \frac{5}{6}\right]; \quad \alpha=1/2\\
\frac{1}{3}\frac{\alpha^2}{2\alpha - 1}\frac{\tau}{T}; \quad \alpha>1/2,
\end{cases}
\end{equation}
where $C\left(\alpha\right) = \frac{\left(1-\alpha\right)\left(2-\alpha\right)\mathbb{B}\left(\alpha + 2, 1 - 2\alpha\right) - 2\left(\alpha^2 + \alpha - 1\right)}{2\left(\alpha+1\right)^2 \left(\alpha + 1\right)}$.
\begin{fact}
Let $\left\{B_{\alpha}\left(t\right)\right\}_{t\geq0}$ be the SBM. Let $\tau_b = \inf\left\{t \geq 0: B_{\alpha}\left(t\right) = b\right\}$ be the time of first hitting a point $b>0$. Then, the pdf of $\tau_b$ is given by \cite{Bodrova2019}
\begin{equation}\label{ht}
    f_{\tau_b}\left(t\right) = \frac{\alpha b}{\sqrt{2\pi}} e^{-\frac{b^2}{2t^{\alpha}}} t^{-1 - \frac{\alpha}{2}}, \quad t \geq 0.
\end{equation}
\end{fact}
\begin{theorem}
    Let $\left\{B_{\alpha}\left(t\right)\right\}_{t\geq0}$ be the SBM. Then, SBM is a martingale with respect to its natural filtration. Its quadratic variation $\left [ B_{\alpha}\right]_t = t^{\alpha}$. The stochastic exponential of SBM is also a martingale with respect to its natural filtration.
    \begin{proof}
        We start with the proof that SBM is a martingale. The first moment of absolute value of SBM is, obviously, finite. For $s<t$ we have:
        \begin{equation}
            \mathbb{E}\left[\giventhat{B_{\alpha}\left(t\right)}{\mathcal{F}_s}\right] = \mathbb{E}\left[\giventhat{B_{\alpha}\left(t\right) - B_{\alpha}\left(s\right) + B_{\alpha}\left(s\right)}{\mathcal{F}_s}\right] = \mathbb{E}\left[B_{\alpha}\left(t\right) - B_{\alpha}\left(s\right)\right] + B_{\alpha}\left(s\right) = B_{\alpha}\left(s\right).
        \end{equation}
        Next, using the SDE representation \eqref{sdesbm}, we calculate the quadratic variation of the SBM. Namely, we have:
        \begin{equation}
            \left [ B_{\alpha}\right]_t = \int_0^t \left(\sqrt{\alpha} s^{\frac{\alpha - 1}{2}}\right)^2 ds = \int_0^t \alpha s^{\alpha - 1} ds = t^{\alpha}.
        \end{equation}
        Now, we will show that the stochastic exponential of SBM is a martingale. Indeed, we have:
        \begin{equation}
            \mathbb{E}\left[\giventhat{e^{B_{\alpha}\left(t\right)}}{\mathcal{F}_s}\right] = \mathbb{E}\left[\giventhat{e^{B_{\alpha}\left(t\right) - B_{\alpha}\left(s\right)} e^{B_{\alpha}\left(s\right)}}{\mathcal{F}_s}\right] =e^{B_{\alpha}\left(s\right)} \mathbb{E}\left[e^{B_{\alpha}\left(t\right) - B_{\alpha}\left(s\right)}\right] = e^{\frac{t^{\alpha} - s^{\alpha}}{2}} e^{B_{\alpha}\left(s\right)}.
        \end{equation}
        The last equality is a straightforward consequence of the Gaussianity. Finally, we have the following:
        \begin{equation}
            \mathbb{E}\left[\giventhat{e^{B_{\alpha}\left(t\right) - \frac{t^{\alpha}}{2}}}{\mathcal{F}_s}\right] = e^{B_{\alpha}\left(s\right) - \frac{s^{\alpha}}{2}}.
        \end{equation}
        As a consequence, it also gives us the integrability of stochastic exponential.
    \end{proof}
\end{theorem}
\section{Scaled Brownian motion with random anomalous diffusion exponent}\label{sbmresec}
\noindent In this section, we define scaled Brownian motion with random anomalous diffusion exponent and investigate its main properties.
\begin{definition}
Let $\mathcal{A}$ be some positive random variable with values from the interval $\left(0, K\right)$ with pdf $f_{\mathcal{A}}\left(\cdot\right)$ and a moment generating function $M_\mathcal{A}\left(\cdot\right)$:
\begin{equation}
    M_\mathcal{A}(s) = \mathbb{E}\left[e^{s\mathcal{A}}\right] =  \int_0^K e^{sa} f_\mathcal{A}\left(a\right)da.
\end{equation}
$K>0$ is finite constant.

\noindent Then, the scaled Brownian motion with random anomalous diffusion exponent (SBMRE) is defined as follows:
\begin{equation}
    B_{\mathcal{A}}\left(t\right) = B\left(t^{\mathcal{A}}\right).
\end{equation}
We assume that $\left\{B\left(t\right)\right\}_{t\geq0}$ and $\mathcal{A}$ are independent. SBMRE can also be represented via the following SDE:

\begin{equation}\label{sdesbmre}
    dB_{\mathcal{A}}\left(t\right) = \sqrt{\mathcal{A}} t^{\frac{\mathcal{A} - 1}{2}}dB\left(t\right), \quad B_{\mathcal{A}}\left(0\right)=0.
\end{equation}
The solution of \eqref{sdesbmre} is given by:
\begin{equation}
    B_{\mathcal{A}}\left(t\right) = \int_0^t \sqrt{\mathcal{A}} s^{\frac{\mathcal{A} - 1}{2}}dB\left(s\right).
\end{equation}
The representation via SDE \eqref{sdesbmre} is a direct consequence of random time change for Brownian motion (for the details see Chapter 8 in \cite{Oeksendal2003}).
\end{definition}
\begin{fact}
    Let $\left\{B_{\mathcal{A}}\left(t\right)\right\}_{t\geq0}$ be the SBMRE. Then, its pdf is given by the following formula:
    \begin{equation}
        p_{B_{\mathcal{A}}}\left(x, t\right) = \int_0^K \frac{1}{\sqrt{2\pi t^a}}exp\left\{\frac{-x^2}{2t^a}\right\} f_\mathcal{A}\left(a\right) da.
    \end{equation}
    \begin{proof}
        It is a direct consequence of the law of total probability.
    \end{proof}
\end{fact}
\begin{fact}\label{moment}
Let $\left\{B_{\mathcal{A}}\left(t\right)\right\}_{t\geq0}$ be the SBMRE. Then, the moment of order q of the absolute value of $\left\{B_{\mathcal{A}}\left(t\right)\right\}_{t\geq0}$ is given by:
\begin{equation}
\label{q-moment}
    \mathbb{E}\left[\left|B_{\mathcal{A}}\left(t\right)\right|^q\right]  = c_q M_\mathcal{A}\left(\frac{1}{2}q\log t\right),
\end{equation}
where $c_q = \frac{2^{q/2} \Gamma\left(\frac{q+1}{2}\right)}{\sqrt{\pi}}$.
\begin{proof}
    \begin{equation}
    \begin{aligned}
        \label{q-moment}
            \mathbb{E}\left[\left|B_{\mathcal{A}}\left(t\right)\right|^q\right] = \mathbb{E}\left[\mathbb{E}\left[\giventhat{\left|B_{\mathcal{A}}\left(t\right)\right|^q} {\mathcal{A}} \right]\right] = \mathbb{E}\left[\mathbb{E}\left[\giventhat{\left|B\left(t^\mathcal{A}\right)\right|^q} {\mathcal{A}} \right]\right]
            = \mathbb{E}\left[\mathbb{E}\left[\giventhat{t^{\frac{1}{2}\mathcal{A}q}\left|B\left(1\right)\right|^q} {\mathcal{A}} \right]\right]
            = \mathbb{E}\left[t^{\frac{1}{2}\mathcal{A}q}Z^q \right] = \\
            = c_q \mathbb{E}\left[e^{\frac{1}{2}q\mathcal{A}\log t}\right] = c_q M_\mathcal{A}\left(\frac{1}{2}q\log t\right),
    \end{aligned}
    \end{equation}
    where $Z = \left|B\left(1\right)\right|^q$ and $\Gamma\left(\cdot\right)$ is the gamma function.
    \end{proof}
\end{fact}
\begin{corollary}\label{msd}
The second moment of $\left\{B_{\mathcal{A}}\left(t\right)\right\}_{t\geq0}$ is given by:
\begin{equation}
    \mathbb{E}\left[B_{\mathcal{A}}^2\left(t\right)\right] = M_\mathcal{A}\left(\log t\right).
\end{equation}
\end{corollary}
\begin{fact}
Let $\left\{B_{\mathcal{A}}\left(t\right)\right\}_{t\geq0}$ be the SBMRE. Then, the autocovariance function of SBMRE is as follows:
\begin{equation}    
    Cov\left(B_{\mathcal{A}}\left(s\right), B_{\mathcal{A}}\left(t\right)\right) = M_\mathcal{A}\left(\log\left(\min\left\{s,  t\right\}\right)\right).
\end{equation}
\begin{proof}
    \begin{equation}
    \begin{aligned}
        Cov\left(B_{\mathcal{A}}\left(s\right),B_{\mathcal{A}}\left(t\right)\right) = \mathbb{E}\left[B_{\mathcal{A}}\left(s\right)B_{\mathcal{A}}\left(t\right)\right] = \mathbb{E}\left[\mathbb{E}\left[B_{\mathcal{A}}\left(s\right)B_{\mathcal{A}}\left(t\right) \giventhat \mathcal{A}\right]\right] = \mathbb{E}\left[\min\left\{s^{\mathcal{A}}, t^{\mathcal{A}}\right\}\right] = \\
        = \mathbb{E}\left[\min\left\{s, t\right\}^{\mathcal{A}}\right] =    \mathbb{E}\left[\exp\left\{\mathcal{A}\log\left(\min\left\{s, t\right\}\right\}\right)\right] = M_\mathcal{A}\left(\log\left(\min\left\{s, t\right\}\right)\right).
    \end{aligned}
    \end{equation}
\end{proof}
\end{fact}
\begin{fact}
   Let $\left\{B_{\mathcal{A}}\left(t\right)\right\}_{t\geq0}$ be the SBMRE. Then, the expected value of TAMSD for SBMRE  is given by:
   \begin{equation}\label{tamsdsbmre}
       \mathbb{E}\left[\delta_\mathcal{A}\left(\tau\right)\right] = \frac{1}{T - \tau} \int_0^{T-\tau} \left(M_\mathcal{A}\left(\log \left( t + \tau\right)\right) - M_\mathcal{A}\left(\log t\right)\right) dt.
   \end{equation}
    \begin{proof}
    \begin{equation}
    \begin{aligned}
    \mathbb{E}\left[\delta_\mathcal{A}\left(\tau\right)\right] = \mathbb{E}\left[\frac{1}{T - \tau} \int_0^{T-\tau} \left(B_{\mathcal{A}}\left(t + \tau\right) - B_{\mathcal{A}}\left(t\right)\right)^2 dt\right] 
    = \mathbb{E}\left[\frac{1}{T - \tau} \int_0^{T-\tau} (B^2_{\mathcal{A}}\left(t + \tau\right) - 2B_{\mathcal{A}}\left(t + \tau\right)B_{\mathcal{A}}\left(t\right) + B^2_{\mathcal{A}}\left(t\right) dt\right] = \\
    = \frac{1}{T - \tau} \int_0^{T-\tau} \left(M_\mathcal{A}\left(\log \left( t + \tau\right)\right) - M_\mathcal{A}\left(\log t\right)\right) dt.
    \end{aligned}
    \end{equation}
The change in the order of integration is made according to the Fubini's theorem \cite{Penot2016}.
    \end{proof}
\end{fact}
\begin{fact}\label{fa45}
    Let $\left\{B_{\mathcal{A}}\left(t\right)\right\}_{t\geq0}$ be the SBMRE. Then, the ergodicity breaking parameter of SBMRE is given by:
    \begin{equation}\label{EB_A}
        EB_\mathcal{A}\left(\tau\right) = \frac{\mathcal{N}_\mathcal{A}\left(\tau\right)}{\mathcal{D}_\mathcal{A}\left(\tau\right)},
    \end{equation}
    where
    \begin{equation}\label{eq35}
\begin{aligned}
    \mathcal{N}_\mathcal{A}\left(\tau\right) = Var\left(\delta_{\mathcal{A}}\left(\tau\right)\right)=\int_0^K \frac{4 \tau^{2a + 2}}{\left(T- \tau\right)^2}\left[\frac{\left(T/\tau - 1\right)^{2a + 1}}{2a + 1} + \frac{\left(3a+1\right)\left(T/\tau - 1\right)^{2a + 2}}{2\left(a + 1\right)^2 \left(2a + 1\right)} - \frac{2\left(T/\tau\right)^{a + 1}\left(T/\tau - 1\right)^{a + 1}}{\left(a + 1\right)^2}\right. + \\
        + \left.\frac{\left(T/\tau\right)^{2a + 2}}{2\left(a + 1\right) \left(2a + 1\right)} - \frac{\left(2a^2 + a+1\right)}{2\left(a + 1\right)^2 \left(2a + 1\right)} + \frac{2}{a +1} \int_0^{T/\tau - 1} x^{a + 1}\left(x+1\right)^{a}dx\right] f_\mathcal{A}\left(a\right) da - \\ - \int_0^K \left(\frac{T^{a +1} - \tau^{a + 1} - \left(T - \tau\right)^{a + 1}}{\left(a+1\right)\left(T-\tau\right)}\right)^2 f_{\mathcal{A}}\left(a\right)da 
        + \left(\int_0^K \frac{T^{a +1} - \tau^{a + 1} - \left(T - \tau\right)^{a + 1}}{\left(a+1\right)\left(T-\tau\right)} f_{\mathcal{A}}\left(a\right)da\right)^2
\end{aligned}
\end{equation}
and
        \begin{equation}
            \mathcal{D}_\mathcal{A}\left(\tau\right) =\mathbb{E}\left[\delta_{\mathcal{A}} \left(\tau\right)\right]^2= \frac{1}{(T - \tau)^2} \int_0^{T-\tau}  \int_0^{T-\tau} \left(M_\mathcal{A}\left(\log \left( t_1 + \tau\right)\right) - M_\mathcal{A}\left(\log t_1\right)\right)\left(M_\mathcal{A}\left(\log \left( t_2 + \tau\right)\right) - M_\mathcal{A}\left(\log t_2\right)\right) dt_1 dt_2.
        \end{equation}
    \begin{proof}
    We can rewrite the EB parameter in the following form using the law of total variance \cite{weiss2005course}:
    \begin{equation}\label{for35}
    EB_{\mathcal{A}}\left(\tau\right) = \frac{Var\left(\delta_{\mathcal{A}}\left(\tau\right)\right)}{\mathbb{E}\left[\delta_{\mathcal{A}} \left(\tau\right)\right]^2}
    = \frac{\mathbb{E}\left[Var\left(\giventhat{\delta_{\mathcal{A}}\left(\tau\right)}{\mathcal{A}}\right)\right] + Var\left(\mathbb{E}\left[\giventhat{\delta_{\mathcal{A}}\left(\tau\right)}{\mathcal{A}}\right]\right)}{\mathbb{E}\left[\delta_{\mathcal{A}} \left(\tau\right)\right]^2}.
\end{equation}
To calculate the first term in the numerator we use \eqref{sbmnumerator} and the law of total expectation,
\begin{equation}\label{eq38}
    \begin{aligned}
        \mathbb{E}\left[Var\left(\giventhat{\delta_{\mathcal{A}}\left(\tau\right)}{\mathcal{A}}\right)\right] = \int_0^K \frac{4 \tau^{2a + 2}}{\left(T- \tau\right)^2}\left[\frac{\left(T/\tau - 1\right)^{2a + 1}}{2a + 1} + \frac{\left(3a+1\right)\left(T/\tau - 1\right)^{2a + 2}}{2\left(a + 1\right)^2 \left(2a + 1\right)} - \frac{2\left(T/\tau\right)^{a + 1}\left(T/\tau - 1\right)^{a + 1}}{\left(a + 1\right)^2}\right. + \\
        + \left.\frac{\left(T/\tau\right)^{2a + 2}}{2\left(a + 1\right) \left(2a + 1\right)} - \frac{\left(2a^2 + a+1\right)}{2\left(a + 1\right)^2 \left(2a + 1\right)} + \frac{2}{a +1} \int_0^{T/\tau - 1} x^{a + 1}\left(x+1\right)^{a}dx\right] f_\mathcal{A}\left(a\right) da.
    \end{aligned}
\end{equation}
To calculate the second term in the numerator we use \eqref{tamsdsbm} and the law of total expectation,
\begin{equation}
    \begin{aligned}
        Var\left(\mathbb{E}\left[\giventhat{\delta_{\mathcal{A}}\left(\tau\right)}{\mathcal{A}}\right]\right) = \mathbb{E}\left[\mathbb{E}\left[\giventhat{\delta_{\mathcal{A}}\left(\tau\right)}{\mathcal{A}}\right]^2\right] - \mathbb{E}\left[\mathbb{E}\left[\giventhat{\delta_{\mathcal{A}}\left(\tau\right)}{\mathcal{A}}\right]\right]^2 = \int_0^K \left(\frac{T^{a +1} - \tau^{a + 1} - \left(T - \tau\right)^{a + 1}}{\left(a+1\right)\left(T-\tau\right)}\right)^2 f_{\mathcal{A}}\left(a\right)da - \\
        - \left(\int_0^K \frac{T^{a +1} - \tau^{a + 1} - \left(T - \tau\right)^{a + 1}}{\left(a+1\right)\left(T-\tau\right)} f_{\mathcal{A}}\left(a\right)da\right)^2.
    \end{aligned}
\end{equation}
Finally, we obtain \eqref{eq35}.

        \noindent By taking a square of \eqref{tamsdsbmre} and changing the square of the integral into a double integral, we get        \begin{equation}
            \mathcal{D}_\mathcal{A}\left(\tau\right) = \frac{1}{(T - \tau)^2} \int_0^{T-\tau}  \int_0^{T-\tau} \left(M_\mathcal{A}\left(\log \left( t_1 + \tau\right)\right) - M_\mathcal{A}\left(\log t_1\right)\right)\left(M_\mathcal{A}\left(\log \left( t_2 + \tau\right)\right) - M_\mathcal{A}\left(\log t_2\right)\right) dt_1 dt_2.
        \end{equation}
    \end{proof}
\end{fact}
\noindent We observe, that if $\mathcal{A}$ is constant, then \eqref{EB_A} reduces to \eqref{ebsbm}.

\begin{fact}
  Let $\left\{B_{\mathcal{A}}\left(t\right)\right\}_{t\geq0}$ be the SBMRE. Then, the distribution of the first hitting time in barrier b for SBMRE $\left\{B_{\mathcal{A}}\left(t\right)\right\}_{t\geq0}$  has the following pdf:
    \begin{equation}
        f_{\tau_b}\left(t\right) = \int_0^K \frac{a b}{\sqrt{2\pi}} e^{-\frac{b^2}{2t^{a}}} t^{-1 - \frac{a}{2}} f_{\mathcal{A}}\left(a\right) da.
    \end{equation}
    \begin{proof}
        We use \eqref{ht} and the law of total probability.
    \end{proof}
\end{fact}
\noindent Next we discuss the martingale properties or SBMRE. In the proof, we follow the same steps as in the proof of Theorem $2.1$ in \cite{Magdziarz2010}. First, we define the filtration $\left\{\mathcal{F}_{\tau}\right\}_{\tau \geq 0}$ given by the following formula:
    \begin{equation}\label{filtration}
        \mathcal{F}_{\tau} = \bigcap_{u>\tau} \left\{ \sigma\left( B\left(y\right) \; : \; 0\leq y\leq u \right) \vee \sigma\left( y^{\mathcal{A}} : \;  y\geq 0 \right) \right \}.
    \end{equation}
    $\mathcal{F}_2 \vee \mathcal{F}_2$ denotes the $\sigma-$ algebra generated by the union of the $\sigma-$ algebras $\mathcal{F}_1, \; \mathcal{F}_2$. 
\begin{theorem}
    Let $\left\{B_{\mathcal{A}}\left(t\right)\right\}_{t\geq0}$ be the SBMRE. Suppose that the density of the random variable $\mathcal{A}$ is continuous. Then, $\left\{B_{\mathcal{A}}\left(t\right)\right\}_{t\geq0}$ is a martingale with respect to filtration $\left\{\mathcal{F}_{t^{\mathcal{A}}}\right\}_{t \geq 0}$, where $\left\{\mathcal{F}_{\tau}\right\}_{\tau \geq 0}$ is given by \eqref{filtration}. The quadratic variation of $\left\{B_{\mathcal{A}}\left(t\right)\right\}_{t\geq0}$ is given by $\left[B_{\mathcal{A}}\right]_t = t^{\mathcal{A}}$. Furthermore, the stochastic exponential of $\left\{B_{\mathcal{A}}\left(t\right)\right\}_{t \geq 0}$ is also $\left\{\mathcal{F}_{t^{\mathcal{A}}}\right\}_{t \geq 0}-$martingale.
    \begin{proof}
    We start with observations that $\left\{\mathcal{F}_{\tau}\right\}$ is right-continuous, $B\left(\tau\right)$ is $\left\{\mathcal{F}_{\tau}\right\}-$martingale, and for every fixed $t_0 > 0$ the random variable $t_0^{\mathcal{A}}$ is the stopping time with respect to $\left\{\mathcal{F}_{t^{\mathcal{A}}}\right\}$. These facts ensure that $\mathcal{G}_t = \left\{\mathcal{F}_{t_0^{\mathcal{A}}}\right\}$ is well defined. Next, we introduce a sequence of $\left\{\mathcal{F}_{\tau}\right\}$-stopping times defined by
    \begin{equation}
        T_n = \inf \left\{\tau>0 : \; \left|B\left(\tau\right)\right| = n\right\}.
    \end{equation}
    Using Doob's optional sampling theorem \cite{Revuz1999} we obtain 
       \begin{equation}\label{doobsampling}
        \mathbb{E}\left[\giventhat{B\left(T_n \wedge t^{\mathcal{A}}\right)}{\mathcal{G}_s}\right] = B\left(T_n \wedge s^{\mathcal{A}}\right) \quad s<t.
    \end{equation}
    Let us note that as $n \to \infty$ the left-hand side and the right-hand side of \eqref{doobsampling} converge  to $\mathbb{E}\left[\giventhat{B\left( t^{\mathcal{A}}\right)}{\mathcal{G}_s}\right]$ and $B\left(s^{\mathcal{A}}\right)$, respectively. Thus, $\left\{B_{\mathcal{A}}\left(t\right)\right\}_{t\geq0}$ is $\left\{\mathcal{G}_t\right\}-$ martingale.
    
    \noindent From the SDE representation of SBMRE \eqref{sdesbmre}  we have that the quadratic variation of $\left\{B_{\mathcal{A}}\left(t\right)\right\}_{t\geq0}$ is given by:
    \begin{equation}
        \left [ B_{\mathcal{A}}\right]_t = \int_0^t \left(\sqrt{\mathcal{A}} s^{\frac{\mathcal{A} - 1}{2}}\right)^2 ds = \int_0^t \mathcal{A} s^{\mathcal{A} - 1} ds = t^{\mathcal{A}}.
    \end{equation}
    Now we show that the stochastic exponential $\left\{Y\left(t\right)\right\}_{t\geq0}$ of SBMRE is also a martingale. Using the Proposition $3.4$ from \cite{Revuz1999} we conclude that $\left\{Y\left(t\right)\right\}_{t\geq0}$ is a local martingale. Hence, it is sufficient to show that $\sup_{0\leq u\leq t} Y\left(u\right)$ is integrable to conclude that $\left\{Y\left(t\right)\right\}_{t\geq0}$ is a martingale. We have the following:
    \begin{equation}
        \sup_{0\leq u\leq t} Y\left(u\right) \leq \sup_{0\leq u\leq t} \exp \left\{\lambda B\left(u^{\mathcal{A}}\right)\right\}.
    \end{equation}
    In addition, $\exp \left\{\lambda B\left(u^{\mathcal{A}}\right)\right\}$ is a positive submartingale. Thus, using Doob's maximal inequality, we get:
    \begin{equation}
        \mathbb{E}\left[\left(\sup_{0\leq u\leq t} \exp \left\{\lambda B\left(u^{\mathcal{A}}\right)\right\}\right)^2\right] \leq 4 \mathbb{E}\left[\exp \left\{2\lambda B\left(u^{\mathcal{A}}\right)\right\}\right].
    \end{equation}
    Now, using the properties of Gaussian random variable and the conditioning on $\sigma\left( y^{\mathcal{A}} : \;  y\geq 0 \right)$, we have:
    \begin{equation}
    \begin{aligned}
        \mathbb{E}\left[\exp \left\{2\lambda B\left(t^{\mathcal{A}}\right)\right\}\right] \leq \mathbb{E}\left[\exp \left\{2\lambda^2 t^{\mathcal{A}}\right\}\right] = \sum_{n=0}^{\infty} \frac{\left(2\lambda^2\right)^n \mathbb{E}\left[t^{n\mathcal{A}}\right]}{n!} = \sum_{n=0}^{\infty} \frac{\left(2\lambda^2\right)^n \mathcal{M}_{\mathcal{A}}\left(n\log\left(t\right)\right)}{n!} = \\ = \sum_{n=0}^{\infty} \frac{\left(2\lambda^2\right)^n \int_0^K e^{an\log\left(t\right)} f_{\mathcal{A}}\left(a\right)da}{n!} = \int_0^K f_{\mathcal{A}}\left(a\right) \sum_{n=0}^{\infty} \frac{\left(2\lambda^2\right)^ne^{an\log\left(t\right)}}{n!} da = \int_0^K e^{2t^a \lambda^2} f_{\mathcal{A}}\left(a\right) da \leq
        e^{2t^c \lambda^2}f_{\mathcal{A}}\left(c\right),
    \end{aligned}
    \end{equation}
    where $c \in \left(0, K \right)$. The last inequality is a straightforward consequence of the mean value theorem and the change in order of integration and summation is performed according to the dominated convergence theorem \cite{Penot2016}.  Finally, $\mathbb{E}\left[\sup_{0\leq u\leq t} Y\left(u\right)\right] < \infty$, thus $\left\{Y\left(t\right)\right\}_{t \geq 0}$ is a martingale.
    \end{proof}
\end{theorem}
\section{Selected distributions of anomalous diffusion exponent}\label{examplessec}
\noindent In this section, we apply the general formulas calculated in the previous section to characterize the SBMRE with two-point and beta distributions of random anomalous diffusion exponent.

\begin{example}\label{ex1}
    First, we consider the case of the simplest distribution of the random variable $\mathcal{A}$, that is, the two-point distribution concentrated at two points $A_1, A_1 \in \left(0, K\right)$, $A_1<A_2$ with pdf given by
    \begin{equation}\label{tppdf}
        f_{\mathcal{A}}(a) = p\delta\left(a - A_1\right) + \left(1-p\right)\delta\left(a-A_2\right),
    \end{equation}
where the weight parameter $p\in \left(0,1\right)$ and $\delta\left(\cdot\right)$ is the Dirac delta function.
The moment generating function of the two-point distribution reads as
\begin{equation}
    M_\mathcal{A}(s) = p e^{A_1 s} + \left(1-p\right)e^{A_2 s} .
\end{equation}
The pdf of SBMRE in this case is given by:
\begin{equation}\label{t2}
    p_{B_{\mathcal{A}}}\left(x, t\right) = \frac{p}{\sqrt{2\pi t^{A_1}}}\exp\left\{-\frac{x^2}{2t^{A_1}}\right\} + \frac{1-p}{\sqrt{2\pi t^{A_2}}}\exp\left\{-\frac{x^2}{2t^{A_2}}\right\}, \quad x \in \mathbb{R}.
\end{equation}
Obviously, \eqref{t2} is a mixture of pdfs of two Gaussian distributed random variables with zero mean and variances equal to $t^{A_1}$ and $t^{A_2}$, respectively.

\noindent The second moment of SBMRE has the form
\begin{equation}\label{tpmsd}
    \mathbb{E}\left[B_{\mathcal{A}}^2\left(t\right)\right] = pt^{A_1} + \left(1-p\right)t^{A_2}.
\end{equation}
Notify that for $t \ll 1$ $t^{A_1} \in o\left(t^{A_2}\right)$ and for $t \gg 1$ $t^{A_2} \in o\left(t^{A_1}\right)$. Thus, the asymptotics behave as
\begin{equation}\label{msdtpasympt}
    \mathbb{E}\left[B_{\mathcal{A}}^2\left(t\right)\right] \sim t^{A},
\end{equation}
where $A=A_1$ for $t \rightarrow 0$ and $A=A_2$ for $t \rightarrow \infty$. This effect is called accelerating diffusion \cite{Chechkin2002, Chechkin2008, Chechkin2011}.

\noindent The expectation of TAMSD reads
\begin{equation}\label{etamsdtp}
    \mathbb{E}\left[\delta_\mathcal{A}\left(\tau\right)\right] = \frac{1}{T-\tau}\left(p \frac{T^{A_1 + 1} - \left(T-\tau\right)^{A_1 + 1} - \tau^{A_1 + 1}}{A_1 +  1} + \left(1-p\right)\frac{T^{A_2 + 1} - \left(T-\tau\right)^{A_2 + 1} - \tau^{A_2 + 1}}{A_2 +  1}\right),
\end{equation}
Taking the first three components of Taylor expansion of $(T-\tau)^{A_1+1}$ and $(T-\tau)^{A_2+1}$ for $\tau/T<<1$ we obtain:
\begin{equation}\label{etamsdtpasympt_1}
 \mathbb{E}\left[\delta_\mathcal{A}\left(\tau\right)\right] \sim p \left(\tau T^{A_1 - 1} + \frac{\tau^{A_1 + 1}}{\left(A_1 + 1\right)T} -  A_1 \tau^2 T^{A_1 - 2}\right) + \left(1 - p \right) \left(\tau T^{A_2 - 1} + \frac{\tau^{A_2 + 1}}{\left(A_2 + 1\right)T} - A_2 \tau^2 T^{A_2 - 2} \right).
\end{equation}
Thus, for the asymptotics at $\tau/T << 1$ (the inequality typical for experiments \cite{metzler2014}) we get:
\begin{equation}\label{etamsdtpasympt}
\mathbb{E}\left[\delta_\mathcal{A}\left(\tau\right)\right]
\sim \frac{\tau}{T} \left(pT^{A_1} + (1-p)T^{A_2}\right),
\end{equation}
which indicates that for the trajectories long enough the $T$-dependence of the expectation value is determined by larger exponent. %

\noindent The EB parameter for two-point diffusion of anomalous diffusion exponent is calculated by the use of  formula (\ref{for35}). The asymptotics for fixed value of $\tau$ and large $T$ are as follows, see \ref{ebasymptoticscalc} for the details,

\begin{equation}\label{twopointebasymptsimple}
    EB_{\mathcal{A}}\left(\tau\right)\sim
\begin{cases}
\frac{pC\left(A_1\right) \tau^{2A_1}+ \left(1-p\right)C\left(A_2\right) \tau^{2A_2} + p\left(1-p\right)\left(T^{A_1} - T^{A_2}\right)^2}{\left(pT^{A_1} + (1-p)T^{A_2}\right)^2}; \quad A_1,A_2<1/2\\\\
\frac{pC\left(A_1\right) \tau^{2A_1} +  \frac{1-p}{12}\tau\left[\log\left(\frac{T}{\tau}\right) + 2\log\left(2\right) - \frac{5}{6}\right] +p\left(1-p\right)\left(T^{A_1} - T^{1/2}\right)^2}{\left(pT^{A_1} + (1-p)T^{1/2}\right)^2}; \quad A_1<1/2, \, A_2=1/2\\\\
\frac{pC\left(A_1\right) \tau^{2A_1} + \frac{\left(1-p\right)A_2^2}{3\left(2A_2 - 1\right)}{\tau}T^{2A_2-1} + p\left(1-p\right)\left(T^{A_1} - T^{A_2}\right)^2}{\left(pT^{A_1} + (1-p)T^{A_2}\right)^2}; \quad A_1<1/2, \, A_2>1/2\\\\
\frac{\frac{p}{12}{\tau}\left[\log\left(\frac{T}{\tau}\right) + 2\log\left(2\right) - \frac{5}{6}\right] + \frac{\left(1-p\right)A_2^2}{3\left(2A_2 - 1\right)}{\tau}T^{2A_2-1} + p\left(1-p\right)\left(T^{1/2} - T^{A_2}\right)^2}{ \left(pT^{1/2} + (1-p)T^{A_2}\right)^2}; \quad A_1=1/2, \, A_2>1/2\\\\
\frac{\frac{pA_1^2}{3\left(2A_1 - 1\right)}{\tau}T^{2A_1-1} + \frac{\left(1-p\right)A_2^2}{3\left(2A_2 - 1\right)}{\tau}T^{2A_2-1} + p\left(1-p\right)\left(T^{A_1} - T^{A_2}\right)^2}{ \left(pT^{A_1} + (1-p)T^{A_2}\right)^2}; \quad A_1>1/2, \, A_2>1/2,
\end{cases}
\end{equation}
where $C\left(\alpha\right)$ is the same as in \eqref{ebsbmsympt}. Let us note that for $p=0$ or $p=1$ and in case $A_1=A_2$, the  formula \eqref{twopointebasymptsimple} reduces to the asymptotics of EB for SBM, see  \eqref{ebsbmsympt}. 
Note the appearance of the last term in the numerator with the prefactor $p(1-p)$, which is absent in the SBM case. 
From formula \eqref{twopointebasymptsimple} we can observe the following universal limit of EB parameter for SBMRE with two-point distribution of anomalous diffusion exponent, 
\begin{equation}
\lim_{T \to \infty} EB_{\mathcal{A}}\left(\tau\right) =
\begin{cases}
    \frac{p}{1-p}, \text{ for } p\neq 0 \text{ and } p\neq 1,\\
    0, \text{ otherwise.}
\end{cases}
\end{equation}

In contrast to the SBM case, in the limit of long trajectories the EB parameter does not depend on the sliding window width $\tau$. This interesting property is checked by Monte Carlo simulations in the next section.

\noindent Further, the first hitting time in the barrier $b$ in the considered case is as follows:
\begin{equation}\label{hittingtp}
    f_{\tau_b}\left(t\right) = \frac{p A_1 b}{\sqrt{2\pi }}\exp\left\{-\frac{b^2}{2t^{A_1}}\right\} t^{-1 - \frac{A_1}{2}} + \frac{(1-p)A_2 b}{\sqrt{2\pi}}\exp\left\{-\frac{b^2}{2t^{A_2}}\right\} t^{-1 - \frac{A_2}{2}}.
\end{equation}
Notify that for $t \gg 1$ we have the following: $\exp\left\{-\frac{b^2}{2t^{A_1}}\right\}, \, \exp\left\{-\frac{b^2}{2t^{A_1}}\right\} \in \mathcal{O}\left(1\right)$ and $t^{-1 - \frac{A_2}{2}} \in o\left(t^{-1 - \frac{A_1}{2}}\right)$. Taking into account these facts,  we obtain the following asymptotic formula for \eqref{hittingtp} in case $t \gg 1$:
\begin{equation}\label{hittingtpasympt}
    f_{\tau_b}\left(t\right) \sim t^{-1 - A_1/2}.
\end{equation}
\end{example}
\begin{example}
The second considered distribution is the beta distribution on the interval $\left[A_1, A_2\right]$, $0<A_1<A_2$ with the pdf
\begin{equation}
    f_{\mathcal{A}}(a) = \frac{\left(a - A_1 \right)^{\gamma - 1} \left(A_2 - a\right)^{\beta - 1}}{\mathbb{B}\left(\gamma, \beta\right)\left(A_2 - A_1\right)^{\gamma + \beta - 1}} \mathbbm{1}_{\left[A_1, A_2\right]},
\end{equation}
where $\gamma, \beta>0$ are beta distribution parameters and $\mathbb{B}\left(\cdot, \cdot \right)$ is beta function. The moment generating function of the beta distribution is given by
\begin{equation}
    M_\mathcal{A}(s) = e^{A_1 s} \prescript{}{1}{F_1 \left(\alpha, \alpha + \beta, s\left(A_2 - A_1\right)\right)},
\end{equation}
where $\prescript{}{1}{F_1 \left(\cdot, \cdot, \cdot\right)}$ is a confluent hypergeometric function (we refer to \cite{Beals2010} for more details).

\noindent The pdf of SBMRE in this case is given by:
\begin{equation}\label{pdf_beta}
    p_{B_{\mathcal{A}}}\left(x, t\right) = \frac{1}{\left(A_2 - A_1 \right)^{\gamma + \beta - 1} \mathbb{B}\left(\gamma, \beta \right)} \int_{A_1}^{A_2} \frac{\left(a - A_1\right)^{\gamma - 1} \left(A_2 - a\right)^{\beta - 1} \exp\left\{- \frac{x^2}{2t^a}\right\}}{\sqrt{2\pi t^a}} da, \quad x \in \mathbb{R}.
\end{equation}
In the case considered one can express the second moment of SBMRE as
\begin{equation}\label{msd_beta}
    \mathbb{E}\left[B_{\mathcal{A}}^2\left(t\right)\right] = t^{A_1} \prescript{}{1}{F_1 \left(\gamma, \gamma + \beta, \left(A_2 - A_1\right) \log\left(t\right)\right)}.
\end{equation}
Using the asymptotics of the hypergeometric function  \cite{Beals2010} we obtain that \eqref{msd_beta} behaves asymptotically as
\begin{equation}\label{msdbetaasymptshort}
\mathbb{E}\left[B_{\mathcal{A}}^2\left(t\right)\right] \sim \frac{\Gamma \left(\gamma + \beta \right)t^{A_1}}{\Gamma \left(\gamma\right)\left(\left(A_2 -A_1\right)\log\left(1/t\right)\right)^{\gamma}}
\end{equation}
for $t \ll 1$ and
\begin{equation}\label{msdbetaasymptlong}
\mathbb{E}\left[B_{\mathcal{A}}^2\left(t\right)\right] \sim \frac{\Gamma \left(\gamma + \beta \right)t^{A_2}}{\Gamma \left(\gamma\right)\left(\left(A_2 -A_1\right)\log\left(t\right)\right)^{\beta}}
\end{equation}
for $t \gg 1$. 

\noindent The expectation of TAMSD for SBMRE with beta distributed anomalous diffusion exponent reads
\begin{equation}
    \mathbb{E}\left[\delta_\mathcal{A}\left(\tau\right)\right] = \frac{1}{T - \tau} \int_0^{T-\tau} \left(e^{A_1 \log\left( t + \tau\right)} \prescript{}{1}{F_1 \left(\gamma, \gamma + \beta, \left(A_2 - A_1\right)\log \left(t + \tau\right)\right)} - e^{A_1 \log t} \prescript{}{1}{F_1 \left(\gamma, \gamma + \beta, \left(A_2 - A_1\right)\log t\right)}\right) dt.
\end{equation}
The EB parameter in this case takes the form
    \begin{equation}
        EB_\mathcal{A}\left(\tau\right) = \frac{\mathcal{N}_\mathcal{A}\left(\tau\right)}{\mathcal{D}_\mathcal{A}\left(\tau\right)},
    \end{equation}
where
\begin{equation}
\begin{aligned}
    \mathcal{N}_\mathcal{A}\left(\tau\right) = \frac{1}{{\mathbb{B}\left(\gamma, \beta\right)\left(A_2 - A_1\right)^{\gamma + \beta - 1}}}\int_{A_1}^{A_2} \frac{4 \tau^{2a + 2}}{\left(T- \tau\right)^2}\left[\frac{\left(T/\tau - 1\right)^{2a + 1}}{2a + 1} + \frac{\left(3a+1\right)\left(T/\tau - 1\right)^{2a + 2}}{2\left(a + 1\right)^2 \left(2a + 1\right)} - \frac{2\left(T/\tau\right)^{a + 1}\left(T/\tau - 1\right)^{a + 1}}{\left(a + 1\right)^2}\right. + \\
        + \left.\frac{\left(T/\tau\right)^{2a + 2}}{2\left(a + 1\right) \left(2a + 1\right)} - \frac{\left(2a^2 + a+1\right)}{2\left(a + 1\right)^2 \left(2a + 1\right)} + \frac{2}{a +1} \int_0^{T/\tau - 1} x^{a + 1}\left(x+1\right)^{a}dx\right]\left(a - A_1 \right)^{\gamma - 1} \left(A_2 - a\right)^{\beta - 1} da + \\ + \frac{1}{{\mathbb{B}\left(\gamma, \beta\right)\left(A_2 - A_1\right)^{\gamma + \beta - 1}}} \int_{A_1}^{A_2} \left(\frac{T^{a +1} - \tau^{a + 1} - \left(T - \tau\right)^{a + 1}}{\left(a+1\right)\left(T-\tau\right)}\right)^2 \left(a - A_1 \right)^{\gamma - 1} \left(A_2 - a\right)^{\beta - 1} da 
        - \\ - \left(\frac{1}{{\mathbb{B}\left(\gamma, \beta\right)\left(A_2 - A_1\right)^{\gamma + \beta - 1}}} \int_{A_1}^{A_2} \frac{T^{a +1} - \tau^{a + 1} - \left(T - \tau\right)^{a + 1}}{\left(a+1\right)\left(T-\tau\right)} \left(a - A_1 \right)^{\gamma - 1} \left(A_2 - a\right)^{\beta - 1}da\right)^2.
\end{aligned}
\end{equation}
and
\begin{equation}
\begin{aligned}
    \mathcal{D}_\mathcal{A}\left(\tau\right) = \frac{1}{\left(T-\tau\right)^2}\int_0^{T-\tau}  \int_0^{T-\tau} \left(\left(t_1 + \tau\right)^{A_1} \prescript{}{1}{F_1 \left(\gamma, \gamma + \beta, \left(A_2 - A_1\right) \log\left(t_1 + \tau\right)\right)} - t_1^{A_1} \prescript{}{1}{F_1 \left(\gamma, \gamma + \beta, \left(A_2 - A_1\right) \log\left(t_1\right)\right)}\right)\\
    \left(\left(t_2 + \tau\right)^{A_1} \prescript{}{1}{F_1 \left(\gamma, \gamma + \beta, \left(A_2 - A_1\right) \log\left(t_2 + \tau\right)\right)} - t_2^{A_1} \prescript{}{1}{F_1 \left(\gamma, \gamma + \beta, \left(A_2 - A_1\right) \log\left(t_2\right)\right)}\right) dt_1 dt_2.
\end{aligned}
\end{equation}
The expression for the first hitting time in barrier $b$ is represented by
\begin{equation}\label{hittingbeta}
    f_{\tau_b}\left(t\right)  = \frac{b}{\sqrt{2\pi}\left(A_2 - A_1 \right)^{\gamma + \beta - 1} \mathbb{B}\left(\gamma, \beta \right)} \int_{A_1}^{A_2} a\left(a - A_1\right)^{\gamma - 1} \left(A_2 - a\right)^{\beta - 1} \exp\left\{- \frac{b^2}{2t^a}\right\} t^{-1 - \frac{a}{2}} da.
\end{equation}
To obtain the asymptotics of \eqref{hittingbeta} for $t \gg 1$  let us note that $\exp\left\{- \frac{b^2}{2t^a}\right\} \in \mathcal{O}\left(1\right)$ as a function of $t$, $a$ is bounded by $A_2$, and $t^{-1 - \frac{a}{2}} \in \mathcal{O}\left(t^{-1 - \frac{A_1}{2}}\right)$ for all $a$. Then we calculate the integral $\int_{A_1}^{A_2}\left(a - A_1\right)^{\gamma - 1} \left(A_2 - a\right)^{\beta - 1}  da$, which is equal to a constant independent on $t$. 
Finally, for \eqref{hittingbeta} we have the following asymptotic formula:
\begin{equation}\label{hittingbetasympt}
    f_{\tau_b}\left(t\right) \sim t^{-1 - A_1/2}.
\end{equation}
\end{example}
\section{Numerical analysis}\label{nasec}

\noindent  In this section, we present the results of numerical analysis demonstrating main properties of the SBMRE (discussed theoretically in the previous section) for two-point distribution and beta distribution of the anomalous diffusion exponent. In addition, we compare
the main characteristics of SBMRE with the corresponding properties of SBM. 

In the numerical analysis we assume $A_1=0.5$ and $A_2=1.5$. These parameters are selected to cover both subdiffusive and superdiffusive cases. Additionally, for comparison purposes, we select the same values of $\alpha$ in the SBM case, namely $\alpha=0.5$ and $\alpha=1.5$. For the two-point distribution of the anomalous diffusive exponent, we consider three $p$ values, namely $p=0.1$, $p=0.5$, and $p=0.9$. In such cases, we observe the dominance of $A_2$ (for $p=0.1$), the balance between $A_1$ and $A_2$ (for $p=0.5$), and the dominance of $A_1$ (for $p=0.9$), respectively. In order to enforce similar behavior for the beta distribution of the anomalous diffusive exponent, we consider three scenarios: $\left(\gamma, \beta\right)=\left(0.7, 0.3\right)$, $\left(\gamma, \beta\right)= \left(0.5, 0.5\right)$, and $\left(\gamma, \beta\right)= \left(0.3, 0.7\right)$.  

\begin{figure}[ht!]
\centering
        \includegraphics[width=0.7\textwidth, height=0.3\textheight]{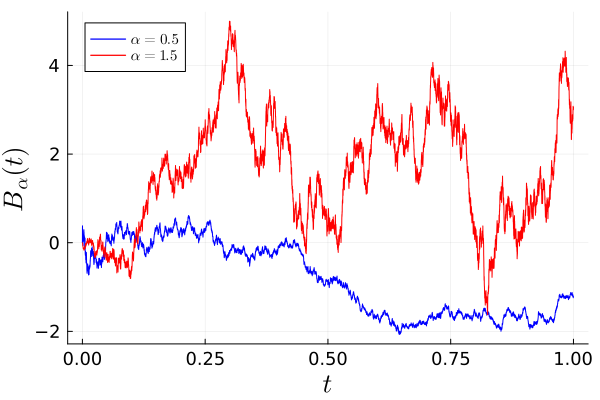}
        \caption{\label{fig:sbmtraj} Exemplary trajectories of SBM  with $\alpha=0.5$ and $\alpha=1.5$. }
\end{figure}

\begin{center}
\begin{figure}[ht!]
       \centering
        \includegraphics[width=1\textwidth, height=0.2\textheight]{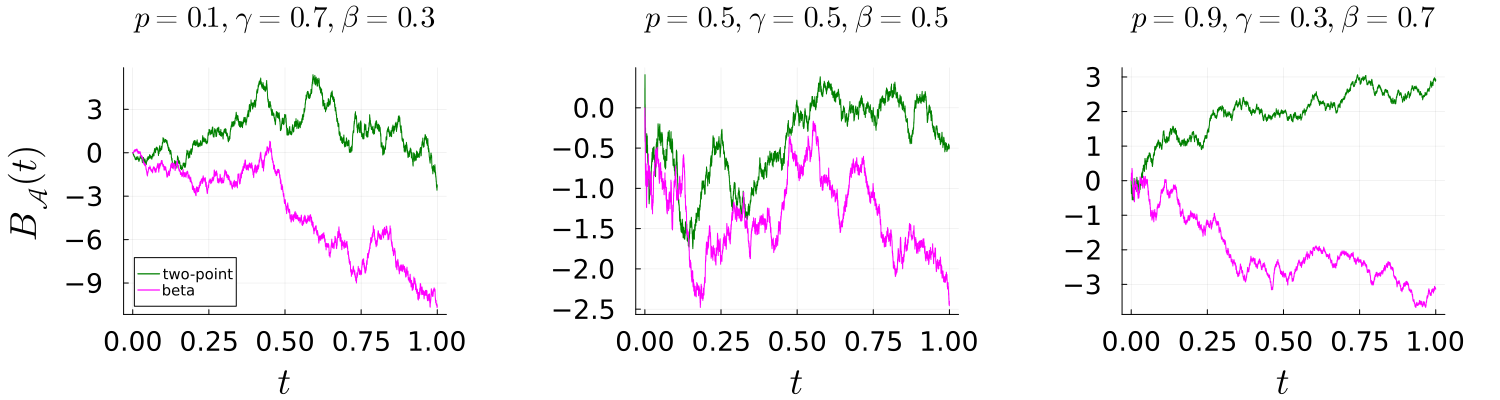}
        \caption{\label{fig:sbmretptraj} Exemplary trajectories of SBMRE with a two-point and beta distribution of the anomalous diffusion exponent with different values of $p$, $\gamma$ and $\beta$; $A_1 = 0.5$, $A_2=1.5$. }
\end{figure}
\end{center}

In Fig. \ref{fig:sbmtraj} we present the exemplary trajectories for SBM for two considered $\alpha$ values,  while in Fig. \ref{fig:sbmretptraj} we demonstrate the trajectories of SBMRE for two-point and beta distributed anomalous diffusive exponent. 
 
\begin{figure}[ht!]
       \centering
        \includegraphics[width=1\textwidth, height=0.2\textheight]{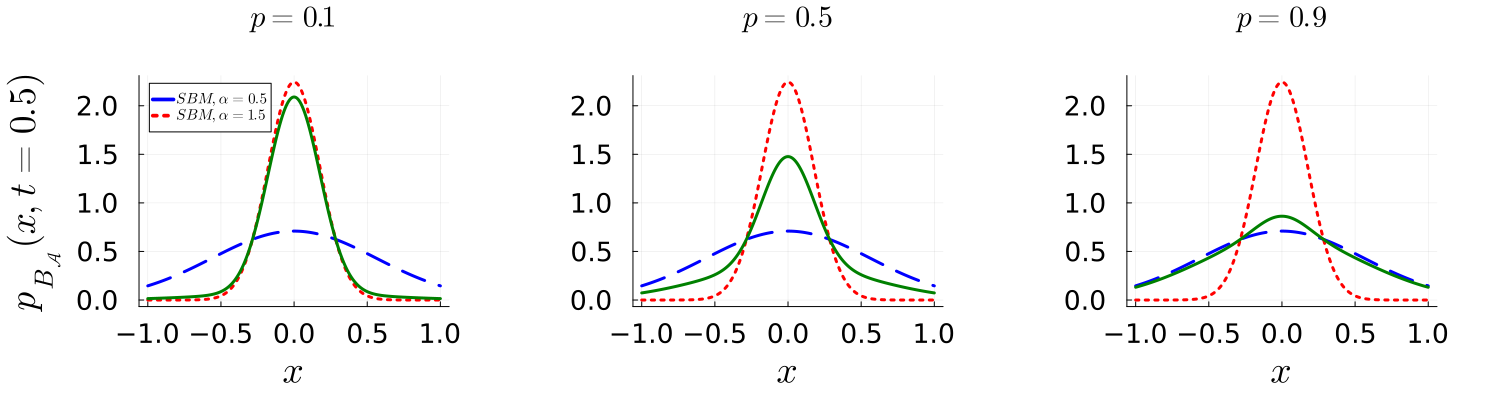} \\\includegraphics[width=1\textwidth, height=0.2\textheight]{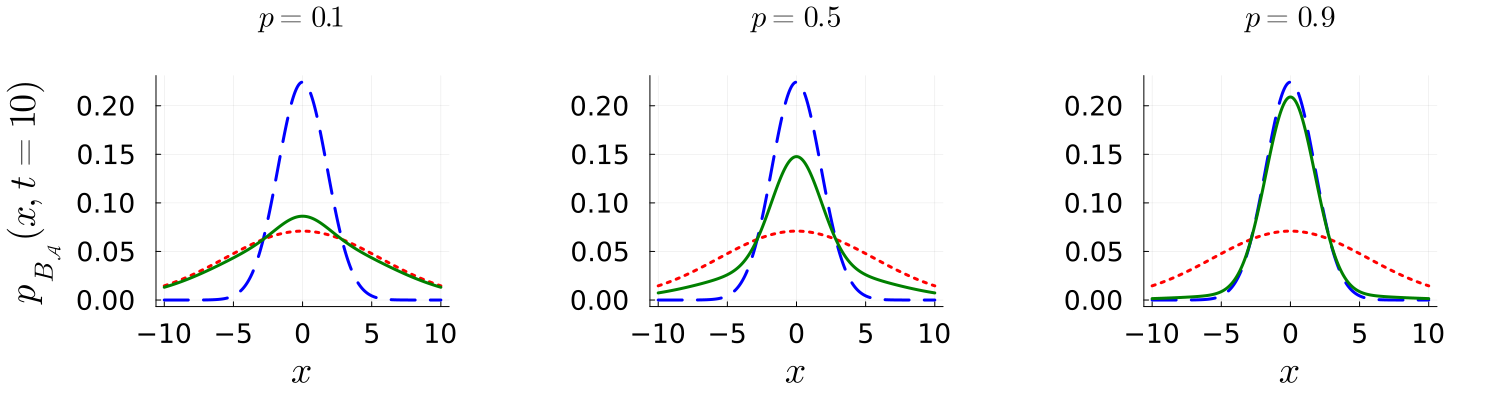}
        \caption{\label{fig:sbmretppdf05}Comparison of pdfs of the position of SBMRE with two-point distribution of the anomalous diffusion exponent and SBM. Top panels: $t=0.5$. Bottom panels: $t=10$. 
        }
\end{figure}

It is evident that basing solely on these plots makes it difficult to distinguish between the specific processes, thus prompting a deeper exploration of their properties. 
\newpage In Fig. \ref{fig:sbmretppdf05}, top panels, we present the pdfs depicting the position of SBMRE  with the two-point distribution of the anomalous diffusion exponent for $t=0.5$, alongside with the pdfs representing the position of subdiffusive and superdiffusive SBM. It is observed that for $p=0.1$ the pdf of SBMRE aligns more closely with that of superdiffusive SBM, and as $p$ increases, it tends towards the subdiffusive counterpart, consistent with intuition. In Fig. \ref{fig:sbmretppdf05}, bottom panels, corresponding pdfs for $t=10$ are presented, where a transition from subdiffusive to superdiffusive cases is apparent. In Fig. \ref{fig_app1} we demonstrate the similar plots for SBMRE with beta distributed anomalous diffusion exponent.

\begin{figure}[ht!]
       \centering
        \includegraphics[width=1\textwidth, height=0.2\textheight]{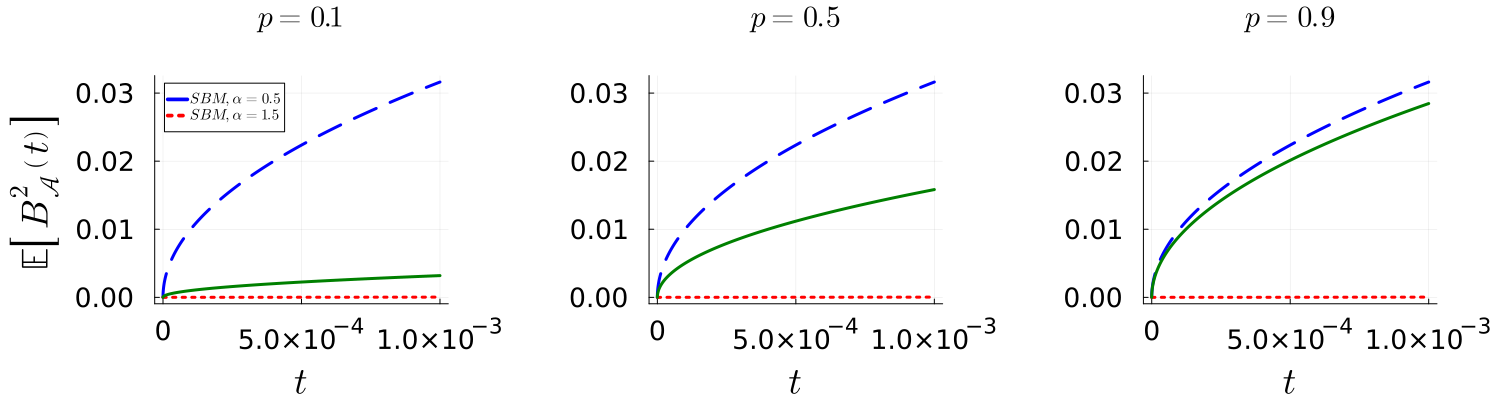}\\
         \includegraphics[width=1\textwidth, height=0.2\textheight]{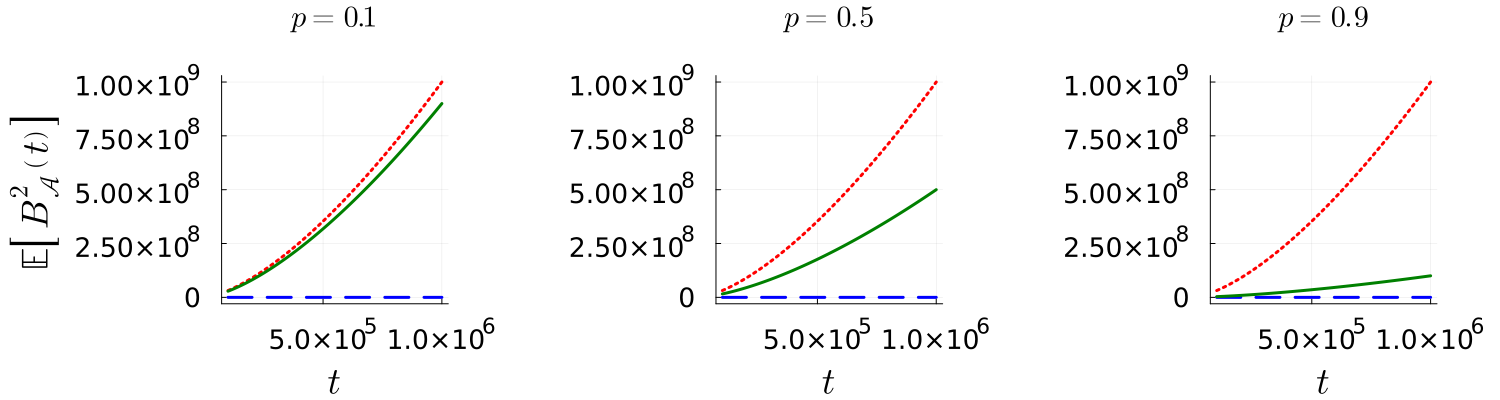}
        \caption{\label{fig:sbmretpmsd}Comparison of the second moments of SBMRE with two-point distribution of the anomalous diffusion exponent and SBM. Top panels: short times Bottom panels: long times.  
        }
\end{figure}

In Fig. \ref{fig:sbmretpmsd}, top panels, we demonstrate short-time behavior of the second moment of SBMRE with a two-point distribution of the anomalous diffusion exponent  exhibiting characteristics akin to both subdiffusive and superdiffusive SBM (depending on $p$ parameter). Notably, we discern that at $p=0.1$, the SBMRE aligns more closely with the superdiffusive scenario of SBM. As the parameter $p$ increases, there is a noticeable shift toward proximity with the subdiffusive scenario. In Fig. \ref{fig:sbmretpmsd}, bottom panels, a similar trend is observed for long time periods, but notably, the transition is from subdiffusive to superdiffusive SBM cases. The second moment for SBMRE with the beta-distributed exponent and SBM is presented in Fig. \ref{app_fig2}. One can observe similar behavior as in case of the two-point distribution. 

In Fig. \ref{fig:sbmretpmsdasympt} (top panels) we depict the asymptotic formula  \eqref{msdtpasympt} of the second moment for SBMRE  for two-point distribution case for short times together with the exact formula presented in \eqref{tpmsd}. Additionally, in Fig. \ref{fig:sbmretpmsdasympt} (bottom panels) we depict the asymptotic formula \eqref{msdtpasympt} of the second moment for SBMRE  for long times together with the exact formula presented in \eqref{tpmsd}. Similar comparison for beta distribution of the anomalous diffusion exponent is presented in Fig. \ref{app_fig3}.   As can be observed, in all cases the theoretical asymptotics perfectly agree with the numerical analysis.

\begin{figure}[ht!]
       \centering
        \includegraphics[width=1\textwidth, height=0.2\textheight]{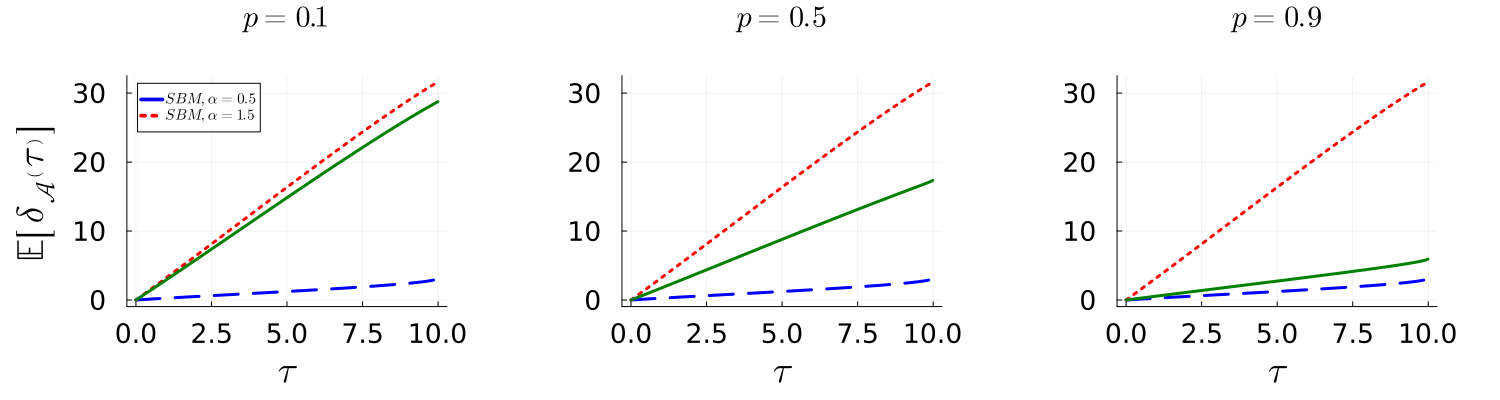}\\
        \caption{\label{fig:sbmretptamsd}Comparison of the expected values of TAMSD of SBMRE with two-point distribution of the anomalous diffusion exponent and SBM. Here $T=10$.
        }
\end{figure}

The expected value of TAMSD for  two-point distribution of the anomalous diffusion exponent is presented in Fig. \ref{fig:sbmretptamsd}. Moreover, in Fig. \ref{fig:sbmretptamsdasympt}, we depict the asymptotic formula \eqref{etamsdtpasympt} of the expected value of the TAMSD for SBMRE with a two-point distribution of the anomalous diffusion exponent together with the exact value presented in \eqref{etamsdtp}. Here we can see the perfect agreement between theoretical formula and values received by numerical analysis. The expected value of TAMSD for the beta distributed anomalous diffusion exponent is visible in Fig. \ref{app_fig4}. The behavior of this statistic for beta distribution is similar to the one observed for two-point distribution for appropriate cases of parameter values.

\begin{figure}[ht!]
       \centering
        \includegraphics[width=1\textwidth, height=0.2\textheight]{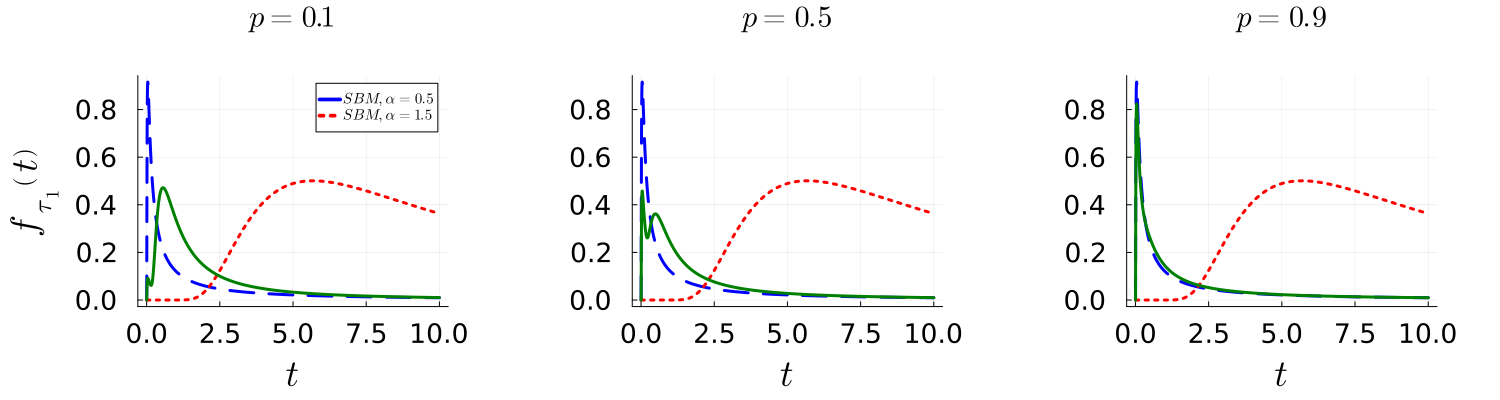}\\
        \caption{\label{fig:sbmretppdfhitting}Comparison of pdf of the hitting time in barrier $b=1$ for SBMRE with two-point distribution of the anomalous diffusion exponent and SBM. 
        }
\end{figure}

Next, in Fig. \ref{fig:sbmretppdfhitting} we plot pdfs of the first hitting time to the barrier $b=1$ for SBMRE  with the two-point distribution and compare the pdf shapes with those for both subdiffusive and superdiffusive SBM. Domination of the lower exponent at long times is clearly visible. For beta distributed anomalous diffusion exponent the same comparison is presented in Fig. \ref{app_fig5}.

\begin{figure}[ht!]
       \centering
         \includegraphics[width=1\textwidth, height=0.2\textheight]{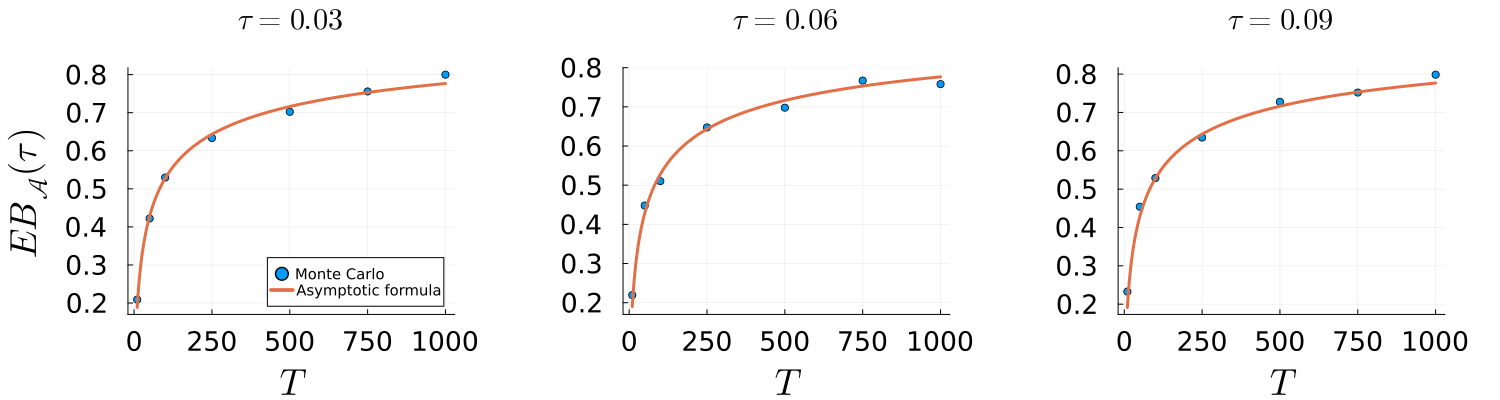}
        \caption{\label{fig::ebmc}EB parameter calculated using 5000 Monte Carlo simulated trajectories of SBMRE with two-point distribution of anomalous diffusion exponent compared with the asymptotic formula \eqref{twopointebasymptsimple} for different values of $\tau$ and $T$; $A_1=0.3$, $A_2=0.7$, $p=0.5$. 
                }
\end{figure}

To demonstrate the interesting behavior of the EB parameter of SBMRE in the case of a two-point distributed anomalous diffusive exponent (discussed in Example \ref{ex1} in Section \ref{examplessec}), in Fig. \ref{fig::ebmc}, we present the EB parameter calculated via Monte Carlo simulations and compare it with the asymptotic formula \eqref{twopointebasymptsimple}. We conducted 5000 Monte Carlo simulations of the SBMRE with a two-point distributed anomalous diffusion exponent with parameters $A_1=0.3$, $A_2=0.7$, and $p=0.5$ (note that the chosen parameters are different from those used in the first part of the presented numerical analyses due to the fact that the case of $A_1=0.5$ seems less interesting in the context of the EB asymptotics).  The simulations were conducted for different values of $\tau$ and $T$. It is noteworthy that there is a very good agreement between the EB obtained from the simulations and its asymptotic behavior at $\tau/T \ll 1$.

\section{Conclusions}
In this paper, we investigate the probabilistic properties of SBMRE, scaled Brownian motion with random anomalous diffusion exponent. This process maintains the characteristics of SBM at the individual trajectory level, albeit possesses randomly varying exponent across the trajectories. SBM is recognized as one of the generic random processes with anomalous diffusion behavior, which finds various applications. The motivation for our study comes from modern biophysical experiments, which clearly demonstrate that standard anomalous diffusion processes with fixed parameters are insufficient to describe the observed phenomena. Consequently, a natural modification of SBM is to incorporate random parameters into the "standard" anomalous diffusion processes, a scenario previously discussed in the context of fractional Brownian motion leading to FBM with a random Hurst exponent \cite{Balcerek2022}, and very recently in doubly stochastic version of continuous time random walk \cite{PhysRevResearch.6.L012033} and L\'evy-walk-like Langevin dynamics with random parameters \cite{10.1063/5.0174613}.

Here, we expand upon the methodology of anomalous diffusion processes with random parameters by analyzing SBMRE. We established the mathematical framework of this process by examining its pdf,  and the q-th absolute moment for the general distribution of the anomalous diffusion exponent. Additionally, we presented the expected value of the TAMSD of SBMRE and discussed its ergodicity breaking property. Furthermore, we explored the properties of the first hitting time pdf and we proven that SBMRE is a martingale, similarly as its stochastic exponential. As special cases of the distributions of the anomalous diffusion exponent, we engage the two-point and beta distributions. In such instances, we discuss the asymptotics of the characteristics of SBMRE. Theoretical findings are validated through numerical analysis, and compared with similar properties of SBM. Similarly to FBMRE, the diffusion of SBMRE is determined by smaller exponent at short times while for sufficiently large times the larger exponent is  dominant, the effect called accelerating diffusion. As for the expected value of TAMSD routinely measured in experiments, the dependence on the trajectory length is determined by larger exponent if the trajectory is long enough. We also observed that the long-time asymptotics of the first hitting time pdf of SBMRE is determined by smaller exponent. Additionally, for the two-point distribution of the anomalous diffusion exponent, the EB parameter for SBMRE does not depend on the sliding window width in the limit of long trajectories, contrasting with the SBM case. A more detailed analysis of this interesting phenomenon and its extension to a more realistic distributions will be provided in a forthcoming publication.  

It is known that the characterization of the anomalous diffusion model from the measurement of individual trajectory is a challenging task that requires new approaches based, in particular, on machine learning or Bayesian inference \cite{nature}. Recently, the latter one was implemented for SBM in \cite{Thapa_2024} and combined with the implementation of Bayesian inference for FBM presented in \cite{Krog2018} to demonstrate model selection and inference of the anomalous diffusion exponent for these two models.  Also, the statistical-based approach for identification of such models was proposed in \cite{Balcerek2021}. It would be of interest to extend these methods for the analysis of anomalous diffusion models with random parameters, such as SBMRE and FBMRE.

\section*{Acknowledgments}
\noindent The work of AW was supported by National Center of Science under Opus Grant 2020/37/B/HS4/00120 “Market risk model identification and validation using novel statistical, probabilistic, and machine learning tools”. 
\bibliographystyle{unsrt} 
\bibliography{bibliography.bib}
\appendix\label{append}
\section{Derivation of asymptotics for EB parameter for two-point distribution of anomalous diffusion exponent.}\label{ebasymptoticscalc}
\noindent Here we derive the asymptotic formulas for EB parameter for two-point distribution of anomalous diffusion exponent \eqref{for35}, when $\tau/T << 1$. Separately, we calculate the asymptotics for the functions given in \eqref{ebtpnumer1}, \eqref{ebtpnumer2}, and \eqref{ebtpdenom}.
Let us note that
\begin{equation}\label{EVar}
    \mathbb{E}\left[Var\left(\giventhat{\delta_{\mathcal{A}}\left(\tau\right)}{\mathcal{A}}\right)\right] = pEB_{A_1}\left(\tau\right)\mathbb{E}\left[\delta_{A_1}\left(\tau\right)\right]^2 + \left(1-p\right)EB_{A_2}\left(\tau\right)\mathbb{E}\left[\delta_{A_2}\left(\tau\right)\right]^2,
\end{equation}
where $EB_{A_1}\left(\tau\right)$, $EB_{A_2}\left(\tau\right)$ are the EB parameters for SBM with $\alpha=A_1$, $\alpha=A_2$ respectively, and $\mathbb{E}\left[\delta_{A_1}\left(\tau\right)\right]^2$, $\mathbb{E}\left[\delta_{A_2}\left(\tau\right)\right]^2$ are squares of expectations of TAMSD for SBM with $\alpha=A_1$, $\alpha=A_2$ respectively.\\
From \eqref{eq38} it follows that
\begin{equation}\label{ebtpnumer1}
    \begin{aligned}
        \mathbb{E}\left[Var\left(\giventhat{\delta_{\mathcal{A}}\left(\tau\right)}{\mathcal{A}}\right)\right] = p\frac{4 \tau^{2A_1 + 2}}{\left(T- \tau\right)^2}\left[\frac{\left(T/\tau - 1\right)^{2A_1 + 1}}{2A_1 + 1} + \frac{\left(3A_1+1\right)\left(T/\tau - 1\right)^{2A_1 + 2}}{2\left(A_1 + 1\right)^2 \left(2A_1 + 1\right)} - \frac{2\left(T/\tau\right)^{A_1 + 1}\left(T/\tau - 1\right)^{A_1 + 1}}{\left(A_1 + 1\right)^2}\right. + \\
        + \left.\frac{\left(T/\tau\right)^{2A_1 + 2}}{2\left(A_1 + 1\right) \left(2A_1 + 1\right)} - \frac{\left(2A_1^2 + A_1+1\right)}{2\left(A_1 + 1\right)^2 \left(2A_1 + 1\right)} + \frac{2}{A_1 +1} \int_0^{T/\tau - 1} x^{A_1 + 1}\left(x+1\right)^{A_1}dx\right] + \\ + \left(1-p\right)\left(\frac{4 \tau^{2A_2 + 2}}{\left(T- \tau\right)^2}\left[\frac{\left(T/\tau - 1\right)^{2A_2 + 1}}{2A_2 + 1} + \frac{\left(3A_2+1\right)\left(T/\tau - 1\right)^{2A_2 + 2}}{2\left(A_2 + 1\right)^2 \left(2A_2 + 1\right)} - \frac{2\left(T/\tau\right)^{A_2 + 1}\left(T/\tau - 1\right)^{A_2 + 1}}{\left(A_2 + 1\right)^2}\right.\right. + \\
        + \left.\left.\frac{\left(T/\tau\right)^{2A_2 + 2}}{2\left(A_2 + 1\right) \left(2A_2 + 1\right)} - \frac{\left(2A_2^2 + A_2+1\right)}{2\left(A_2 + 1\right)^2 \left(2A_2 + 1\right)} + \frac{2}{A_2 +1} \int_0^{T/\tau - 1} x^{A_2 + 1}\left(x+1\right)^{A_2}dx\right]\right).
    \end{aligned}
\end{equation}
Additionally, we have
\begin{equation}\label{ebtpnumer2}
    \begin{aligned}
        Var\left(\mathbb{E}\left[\giventhat{\delta_{\mathcal{A}}\left(\tau\right)}{\mathcal{A}}\right]\right) =\mathbb{E}\left[\mathbb{E}\left[\giventhat{\delta_{\mathcal{A}}\left(\tau\right)}{\mathcal{A}}\right]^2\right] - \mathbb{E}\left[\mathbb{E}\left[\giventhat{\delta_{\mathcal{A}}\left(\tau\right)}{\mathcal{A}}\right]\right]^2  = p\left(\frac{T^{A_1 +1} - \tau^{A_1 + 1} - \left(T - \tau\right)^{A_1 + 1}}{\left(A_1+1\right)\left(T-\tau\right)}\right)^2  + \\ + \left(1-p\right)\left(\left(\frac{T^{A_2 +1} - \tau^{A_2 + 1} - \left(T - \tau\right)^{A_2 + 1}}{\left(A_2+1\right)\left(T-\tau\right)}\right)^2 \right)
        - \left(p\frac{T^{A_1 +1} - \tau^{A_1 + 1} - \left(T - \tau\right)^{A_1 + 1}}{\left(A_1+1\right)\left(T-\tau\right)} + \left(1-p\right)\left(\frac{T^{A_2 +1} - \tau^{A_2 + 1} - \left(T - \tau\right)^{A_2 + 1}}{\left(A_2+1\right)\left(T-\tau\right)}\right)\right)^2.
    \end{aligned}
\end{equation}
Finally,
\begin{equation}\label{ebtpdenom}
          \mathbb{E}\left[\delta_{\mathcal{A}} \left(\tau\right)\right]^2= \left(p\frac{T^{A_1 +1} - \tau^{A_1 + 1} - \left(T - \tau\right)^{A_1 + 1}}{\left(A_1+1\right)\left(T-\tau\right)} + \left(1-p\right)\left(\frac{T^{A_2 +1} - \tau^{A_2 + 1} - \left(T - \tau\right)^{A_2 + 1}}{\left(A_2+1\right)\left(T-\tau\right)}\right)\right)^2.
        \end{equation}
\noindent To obtain \eqref{ebtpnumer1}, \eqref{ebtpnumer2}, \eqref{ebtpdenom} we use formula \eqref{tppdf} for the pdf of two-point distribution. \\

\noindent Let us note that the asymptotics of  $\mathbb{E}\left[Var\left(\giventhat{\delta_{\mathcal{A}}\left(\tau\right)}{\mathcal{A}}\right)\right]$  is equal to the asymptotics of the EB parameter for SBM (see \eqref{ebsbmsympt}) multiplied by the square of the asymptotic of expectation of TAMSD for SBM, see \eqref{tamsdsbmasympt}.  Since the asymptotic of the EB parameter is different for different ranges of parameter $\alpha$ we have to divide our analysis into few cases.
\noindent First, we assume $A_1,A_2<1/2$. In this case \eqref{EVar} takes the form
\begin{equation}\label{ebasympttp1}
    \mathbb{E}\left[Var\left(\giventhat{\delta_{\mathcal{A}}\left(\tau\right)}{\mathcal{A}}\right)\right] \sim pC\left(A_1\right) \frac{\tau^{2A_1+2}}{T^{2}} + \left(1-p\right)C\left(A_2\right) \frac{\tau^{2A_2+2}}{T^{2}},
\end{equation}
\\
Further, for $A_1<1/2$ and $A_2 = 1/2$
\begin{equation}
    \mathbb{E}\left[Var\left(\giventhat{\delta_{\mathcal{A}}\left(\tau\right)}{\mathcal{A}}\right)\right] \sim pC\left(A_1\right) \frac{\tau^{2A_1+2}}{T^{2}} +  \left(1-p\right)\frac{\tau^3}{12T^{3 - 2A_2}}\left[\log\left(\frac{T}{\tau}\right) + 2\log\left(2\right) - \frac{5}{6}\right],
\end{equation}
for $A_1<1/2$ and $A_2>1/2$
\begin{equation}
    \mathbb{E}\left[Var\left(\giventhat{\delta_{\mathcal{A}}\left(\tau\right)}{\mathcal{A}}\right)\right] \sim pC\left(A_1\right) \frac{\tau^{2A_1+2}}{T^{2}} +\left(1-p\right)\frac{A_2^2}{3\left(2A_2 - 1\right)}\frac{\tau^3}{T^{3-2A_2}},
\end{equation}
\noindent and for $A_1=1/2$ and $A_2 > 1/2$
\begin{equation}
    \mathbb{E}\left[Var\left(\giventhat{\delta_{\mathcal{A}}\left(\tau\right)}{\mathcal{A}}\right)\right] \sim  p\frac{\tau^3}{12T^{3 - 2A_1}}\left[\log\left(\frac{T}{\tau}\right) + 2\log\left(2\right) - \frac{5}{6}\right] + \left(1-p\right)\frac{A_2^2}{3\left(2A_2 - 1\right)}\frac{\tau^3}{T^{3-2A_2}} .
\end{equation}
\noindent Finally, for $A_1, A_2 > 1/2$ the following holds:
\begin{equation}
    \mathbb{E}\left[Var\left(\giventhat{\delta_{\mathcal{A}}\left(\tau\right)}{\mathcal{A}}\right)\right] \sim  p\frac{A_2^2}{3\left(2A_2 - 1\right)}\frac{\tau^3}{T^{3-2A_1}} + \left(1-p\right)\frac{A_2^2}{3\left(2A_2 - 1\right)}\frac{\tau^3}{T^{3-2A_2}} .
\end{equation}
For  $Var\left(\mathbb{E}\left[\giventhat{\delta_{\mathcal{A}}\left(\tau\right)}{\mathcal{A}}\right]\right)$ we have the following asymptotics when $\tau/T<<1$:
\begin{equation}
    \begin{aligned}
      Var\left(\mathbb{E}\left[\giventhat{\delta_{\mathcal{A}}\left(\tau\right)}{\mathcal{A}}\right]\right) \sim \left[\left(\frac{\tau}{T}\right)^2 \left(pT^{2A_1} + (1-p)T^{2A_2}\right)\right] - \left(\frac{\tau}{T} \left(pT^{A_1} + (1-p)T^{A_2}\right)\right)^2 = \\ = \left(\frac{\tau}{T}\right)^2\left[ \left(p-p^2\right)\left(T^{2A_1} + T^{2A_2}\right)-2p\left(1-p\right)T^{A_1 + A_2}\right].  
    \end{aligned}
\end{equation}
The above is the straightforward consequence of the asymptotics of expectation of TAMSD for SBM presented in \eqref{tamsdsbm}. Let us note, in case $A_1=A_2=1/2$ the SBMRE reduces to SBM with $\alpha=1/2$ and thus this case is not considered here.
\noindent Finally, for \eqref{ebtpdenom} the following holds, when $\tau/T<<1$ 
\begin{equation}\label{ebasympttplast}
    \mathbb{E}\left[\delta_\mathcal{A}\left(\tau\right)\right]^2 \sim \left(\frac{\tau}{T} \left(pT^{A_1} + (1-p)T^{A_2}\right)\right)^2,
\end{equation}
which is just the square of the asymptoics of expectation of TAMSD for SBMRE with two-point distribution of anomalous diffsuion exponent, see \eqref{etamsdtpasympt}. Taking formulas \eqref{ebasympttp1}-\eqref{ebasympttplast}  we obtain the asymptotics for the EB parameter for two-point distributed anomalous diffusion exponent given by \eqref{twopointebasymptsimple}. 
\\

\section{Additional plots.}
\begin{figure}[h!]
       \centering
        \includegraphics[width=1\textwidth, height=0.2\textheight]{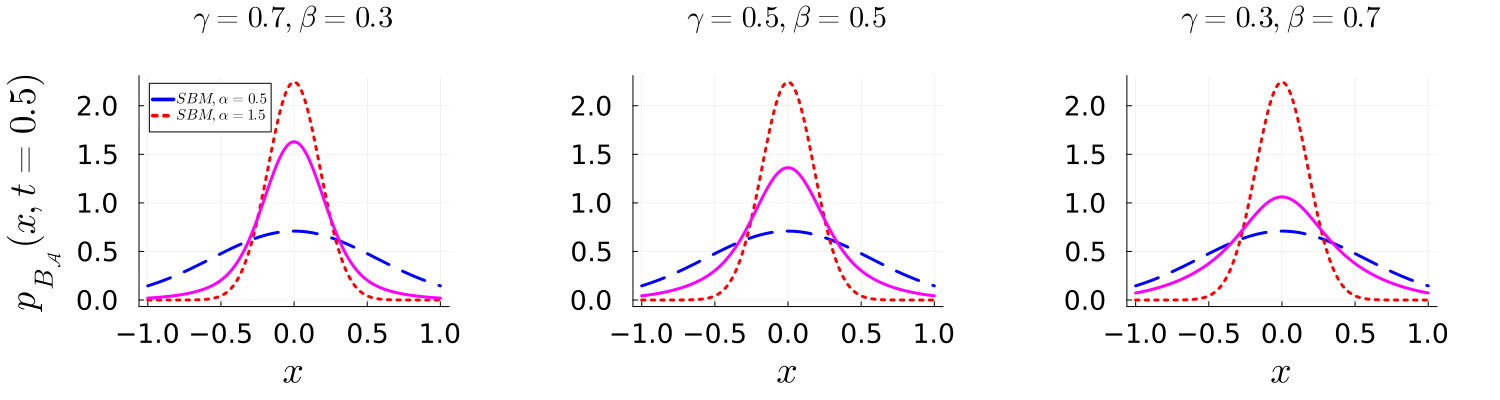}\\  \includegraphics[width=1\textwidth, height=0.2\textheight]{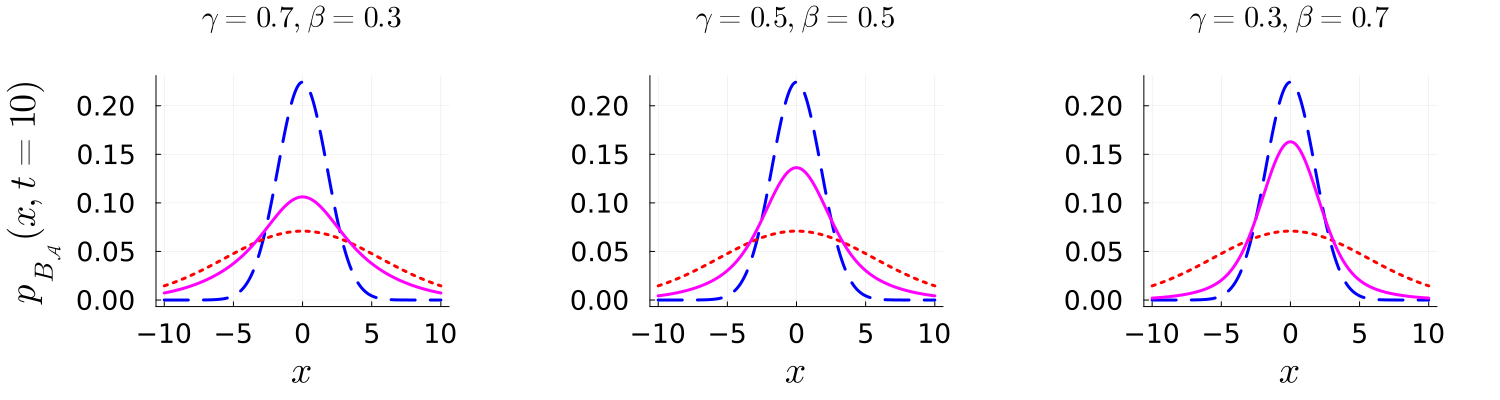}
        \caption{\label{fig:sbmrebetapdf}Comparison of pdfs of the position of SBMRE with beta distribution of the anomalous diffusion exponent and SBM. Top panels: $t=0.5$.  Bottom panels:  $t=10$.  
        }\label{fig_app1}
\end{figure}

\begin{figure}[h!]
       \centering
        \includegraphics[width=1\textwidth, height=0.2\textheight]{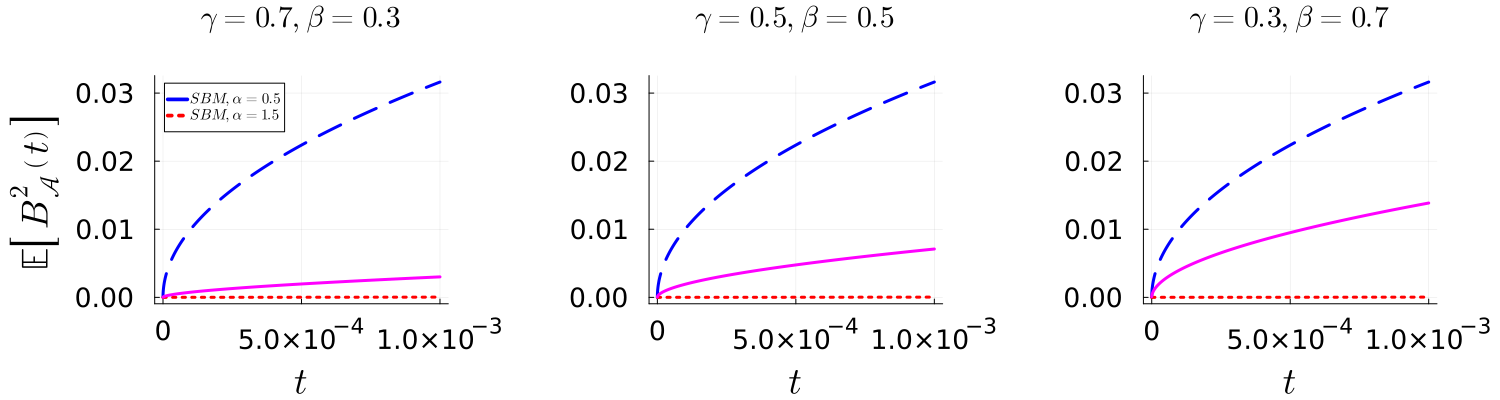}\\  \includegraphics[width=1\textwidth, height=0.2\textheight]{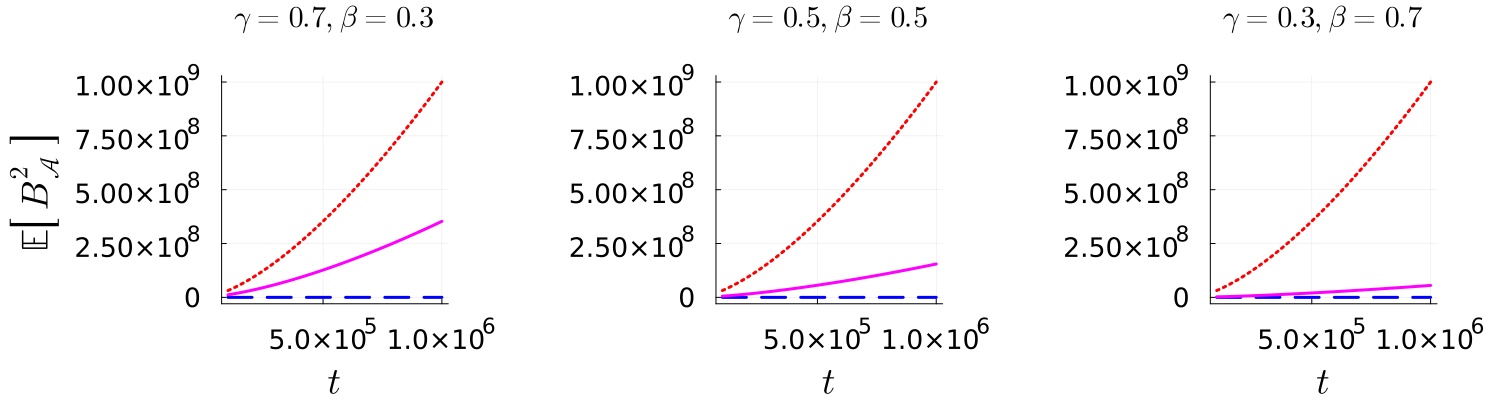}
        \caption{\label{fig:sbmrebetamsdlong} Comparison  of second moments of SBMRE with  beta distribution of the anomalous diffusion exponent and SBM. Top panels: short times. Bottom panels:  long times.  
        }\label{app_fig2}
\end{figure}
\begin{figure}[h!]
       \centering
        \includegraphics[width=1\textwidth, height=0.2\textheight]{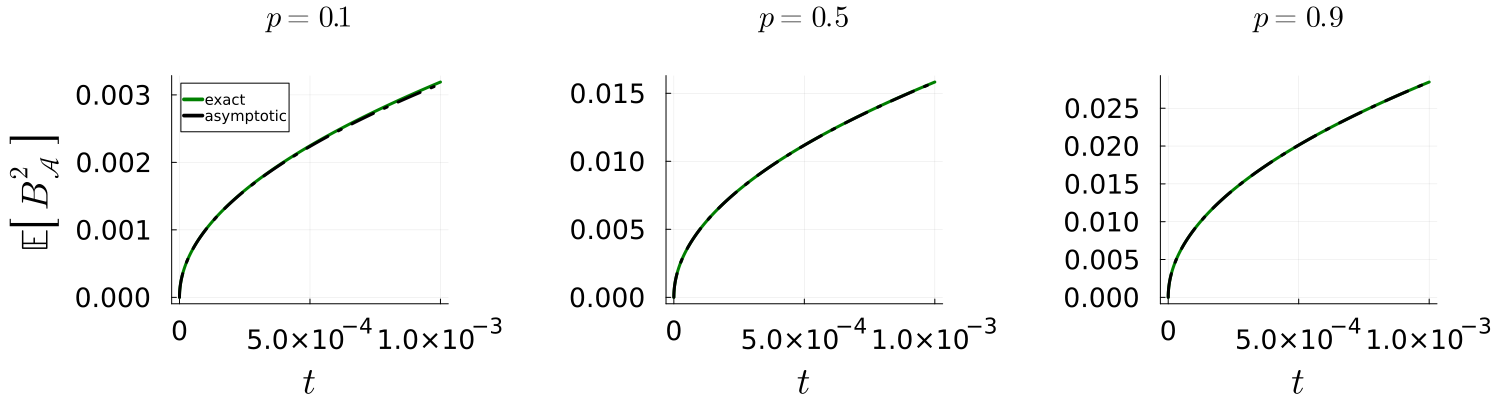}\\  \includegraphics[width=1\textwidth, height=0.2\textheight]{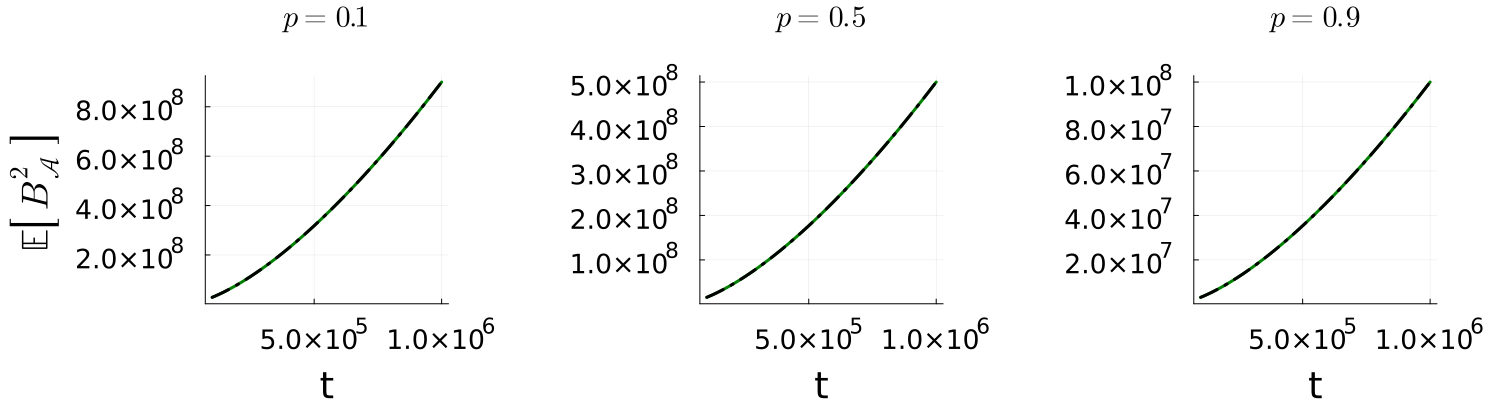}
        \caption{\label{fig:sbmretpmsdasympt} Comparison of second moments of SBMRE with two-point distribution of the anomalous diffusion exponent and their asymptotics (see formula \eqref{msdtpasympt}). Top panels: short times. Bottom panels: long times.  
        }
\end{figure}
\begin{figure}[h!]
       \centering
        \includegraphics[width=1\textwidth, height=0.2\textheight]{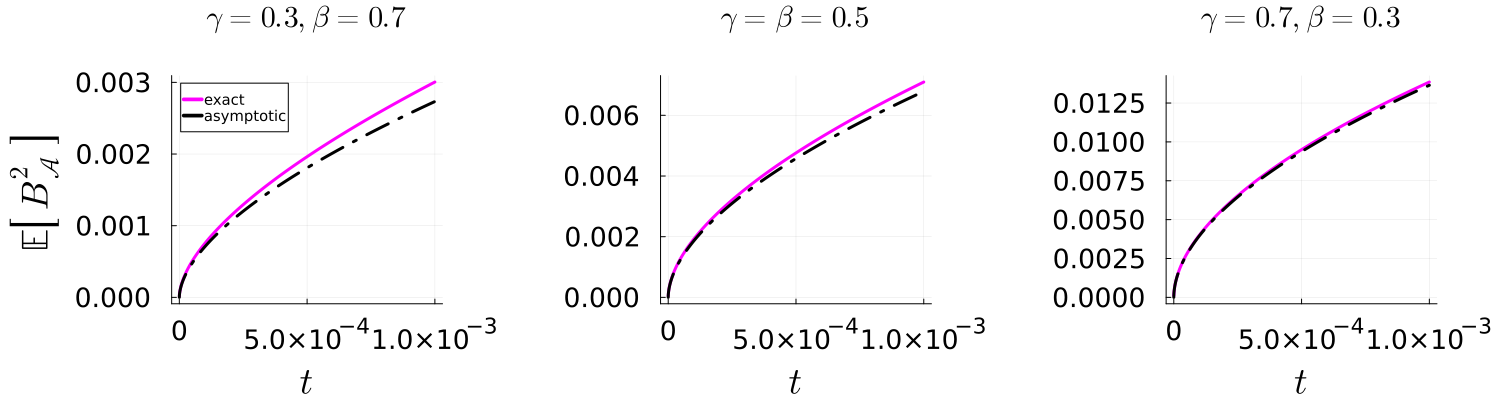}\\  \includegraphics[width=1\textwidth, height=0.2\textheight]{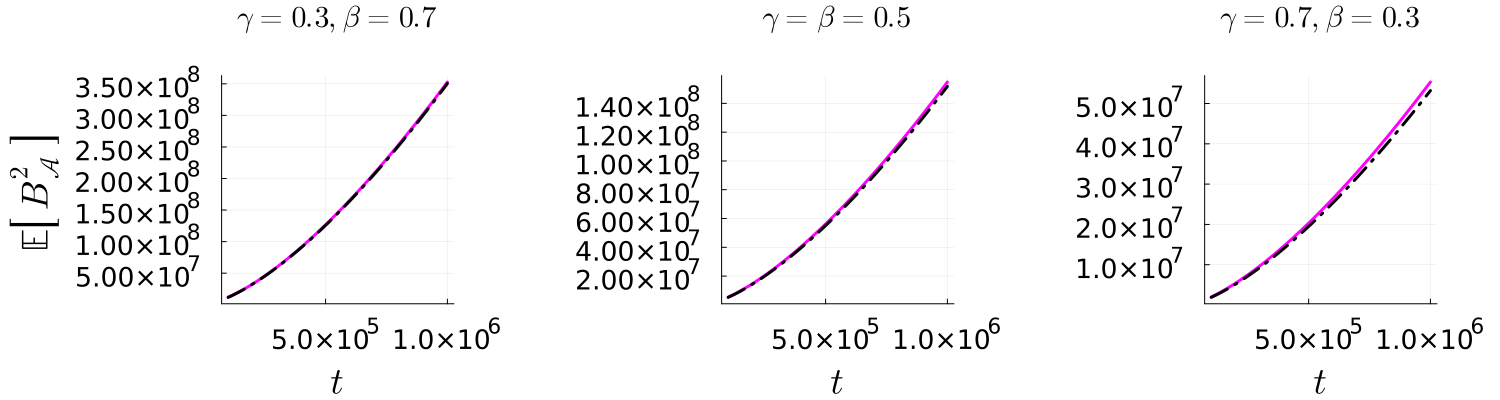}
        \caption{\label{fig:sbmrebetamsdasympt} Comparison of second moments of SBMRE with beta distribution of the anomalous diffusion exponent and their asymptotics (see formulas \eqref{msdbetaasymptshort}, \eqref{msdbetaasymptlong}). Top panels: short times. Bottom panels: long times.
        }\label{app_fig3}
\end{figure}
\begin{figure}[ht!]
       \centering
        \includegraphics[width=1\textwidth]{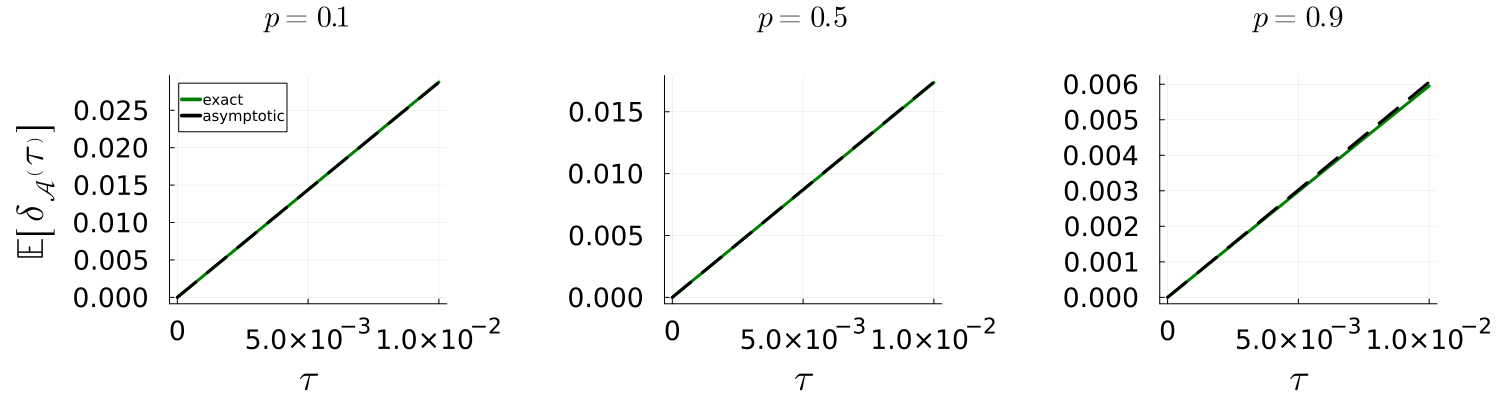}
        \caption{\label{fig:sbmretptamsdasympt}Asymptotic (see formula \eqref{etamsdtpasympt}) of the expectation of the TAMSD of SBMRE with a two-point distribution of the anomalous diffusion exponent. Here $T=10$. 
        }\label{app_fig10}
\end{figure}
\begin{figure}[ht!]
       \centering
          \includegraphics[width=1\textwidth, height=0.2\textheight]{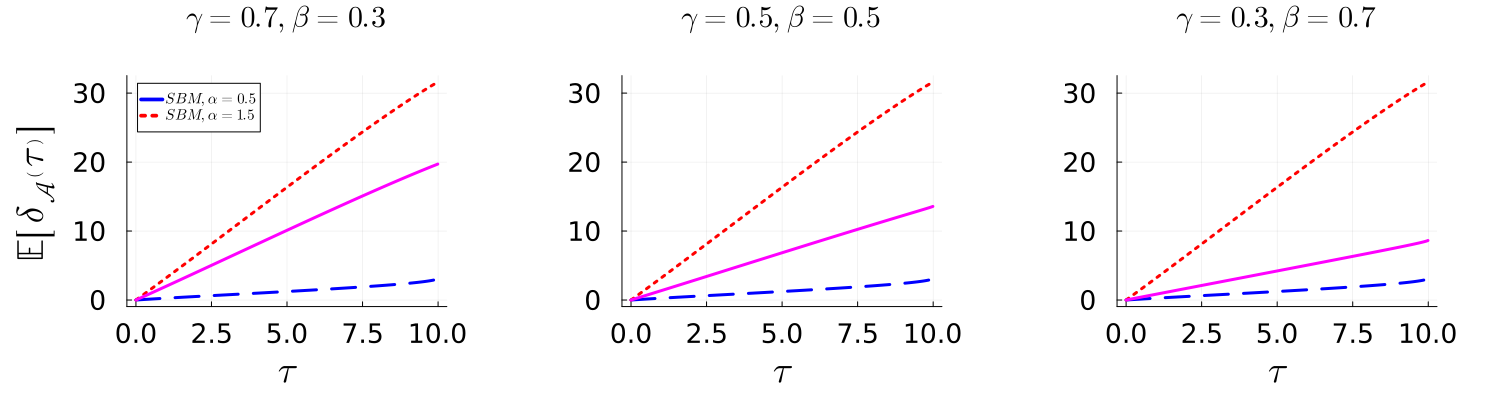}
        \caption{\label{fig:sbmrebetatamsd}Comparison of the expected values of TAMSD of SBMRE with beta distribution of the anomalous diffusion exponent and SBM. Here $T=10$.
        }\label{app_fig4}
\end{figure}

\begin{figure}[ht!]
       \centering
         \includegraphics[width=1\textwidth, height=0.2\textheight]{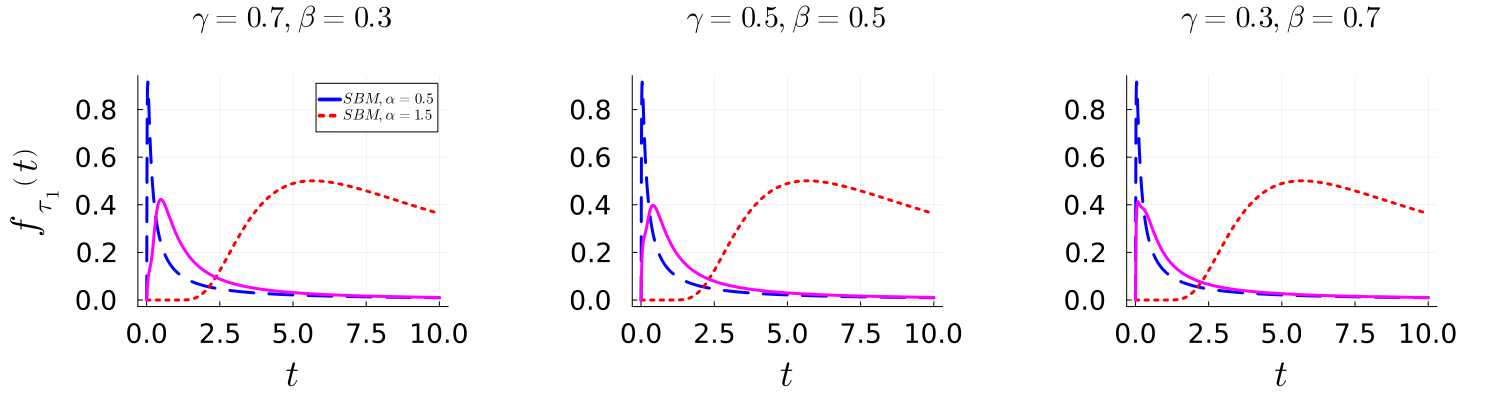}
        \caption{\label{fig:sbmrebetapdfhitting}Comparison of pdf of the hitting time in barrier $b=1$ for SBMRE with beta distribution of the anomalous diffusion exponent and SBM. 
        }\label{app_fig5}
\end{figure}

\end{document}